\documentclass{amsart}

\usepackage{amsmath}
\usepackage{amssymb}
\usepackage{amsfonts}
\usepackage{MnSymbol}
\usepackage{graphics} 
\usepackage{epsfig} 
\usepackage{psfrag}
\usepackage{pdfsync}
\usepackage[usenames]{color}
\usepackage{cancel}

\usepackage[notcite,notref]{showkeys} 

\definecolor{arancio}{rgb}{0.90,0.50,0.20}

\definecolor{blu}{rgb}{0.,0.,1.}

\definecolor{pavone}{rgb}{0.00,0.00,0.63}

\definecolor{malva}{rgb}{0.10,0.50,0.50}

\definecolor{rosso}{rgb}{1.,0.,0.}

\definecolor{geranio}{rgb}{0.90,0.00,0.20}

\definecolor{cerulean}{rgb}
{0.0, 0.48, 0.65}

\newtheorem{theorem}{Theorem}[section]

\newtheorem{lem}[theorem]{Lemma}
\newtheorem{prop}[theorem]{Proposition}

\theoremstyle{definition}
\newtheorem{definition}{Definition}[section]
\newtheorem{remark}{Remark}[section]

\newcommand{\N}{\mathbb{N}}
\newcommand{\R}{\mathbb{R}}
\newcommand{\lla}{\left\langle}
\newcommand{\rra}{\right\rangle}
\date{\today}

\newcommand{\bcl}{\begin{center}}
\newcommand{\ecl}{\end{center}}
\newcommand{\brl}{\begin{right}}
\newcommand{\erl}{\end{right}}
\newcommand{\ben}{\begin{enumerate}}
\newcommand{\barr}{\begin{array}}
\newcommand{\earr}{\end{array}}
\newcommand{\btab}{\begin{tabular}}
\newcommand{\etab}{\end{tabular}}
\newcommand{\bdoc}{\begin{document}}
\newcommand{\edoc}{\end{document}}
\newcommand{\beqy}{\begin{eqnarray}}

\newcommand{\beq}{\begin{equation}}
\newcommand{\beqi}{\begin{eqnarray*}}
\newcommand{\bitem}{\begin{itemize}}
\newcommand{\brem}{\begin{remark}}
\newcommand{\erem}{\end{remark}}
\newcommand{\eitem}{\end{itemize}}
\newcommand{\nln}{\newline}
\newcommand{\newt}{\newtheorem}

\renewcommand{\a }{\alpha }
\renewcommand{\b }{\beta }
\newcommand{\g }{\gamma}
\newcommand{\G }{\Gamma }
\renewcommand{\d }{\delta }

\newcommand{\D }{\Delta }
\newcommand{\e }{\epsilon }
\newcommand{\z }{\zeta }
\renewcommand{\l }{\lambda }
\renewcommand{\L }{\Lambda }
\newcommand{\m }{\mu }
\newcommand{\n }{\tau }
\renewcommand{\r }{\rho }
\newcommand{\s }{\sigma }
\newcommand{\Sig }{\Sigma }
\renewcommand{\t }{\tau }
\newcommand{\var }{\varphi }
\renewcommand{\o }{\omega }
\renewcommand{\O }{\Omega }

\newcommand{\supp}{\text{\rm supp}\,}
\newcommand{\sgn}{\text{\rm sgn}\,}

\newcommand{\red}[1]{{\color{red}{#1}}}
\newcommand{\blue}[1]{{\color{blue}{#1}}}

\title[Measure-valued solutions]
{ A uniqueness criterion for measure-valued solutions of scalar hyperbolic conservation laws}

\author[Bertsch]{Michiel Bertsch}
\address{Dipartimento di Matematica, Universit\`a di Roma  Tor Vergata, 
Via della Ricerca Scientifica, 00133 Roma, Italy \\ and
Istituto per le Applicazioni del Calcolo ``M. Picone", CNR, Roma, Italy} 
\email{bertsch.michiel@gmail.com}
\author[Smarrazzo]{Flavia Smarrazzo}
\address{Universit\`a Campus Bio-Medico di Roma\\ Via Alvaro del Portillo 21, 00128 Roma, Italy}
\thanks{}
\email{flavia.smarrazzo@gmail.com}
\author[Terracina]{Andrea Terracina}
\address{Dipartimento di Matematica ``G. Castelnuovo", Universit\`a Sapienza di Roma\\ P.le A. Moro 5, I-00185 Roma, Italy}
\email{terracina@mat.uniroma1.it}
\author[Tesei]{Alberto Tesei}
\address{Dipartimento di Matematica ``G. Castelnuovo", Universit\`a Sapienza di Roma\\ P.le A. Moro 5, I-00185 Roma, Italy, and 
Istituto per le Applicazioni del Calcolo ``M. Picone", CNR, Roma, Italy}
\email{albertotesei@gmail.com}



\date{\today}

\begin{document}

\bibliographystyle{h-elsevier2}
\begin{abstract} 

We prove existence and uniqueness of Radon measure-valued solutions of the Cauchy problem 
$$
\left\{\begin{array}{ll}
u_t+ \left[\varphi(u)\right]_x=0 
& \quad\mbox{in}\  \R\times (0,T) \\
u=u_0\ge 0 &\quad\mbox{in}\  \R\times \{0\}\,,
\end{array}\right. 
$$ 
where $u_0$ a positive 
 Radon measure whose singular part is a  finite superposition of Dirac masses, and $\varphi\in C^2([0,\infty))$ is bounded. The novelty of the paper is the introduction of a compatibility condition which, combined with standard entropy conditions, guarantees uniqueness. 
\end{abstract}

\maketitle


\section{Introduction}\label{intro}

\subsection{Statement of the problem}

In this paper we consider the Cauchy problem 
$$
\left\{\begin{array}{ll}
u_t+\left[\varphi(u)\right]_x=0 
& \quad\mbox{in}\  \R\times (0,T)=:S \smallskip\\
u=u_0  &\quad\mbox{in}\  \R\times \{0\}\,.
\end{array}\right. \leqno{(P)}
$$  
Here $T>0$, $u_0$ is a positive Radon measure on $\R$ whose singular part $u_{0s}$ (with respect to the Lebesgue measure) is a finite superposition of Dirac masses, and $\varphi$ is a  smooth and bounded function with bounded derivative:
 $$
u_{0s}=\sum_{j=1}^{p} c_j \delta_{x_j}\qquad 
\text{($x_1<x_2<\dots<x_p$;  \ $c_j>0$ for $1\le j\le p$),}
\leqno{(H_0)}
$$
$$
\text{$\varphi\in C^2([0,\infty))\cap W^{1,\infty}(0,\infty)$\,.} \leqno{(H_1)}
$$

It is worth mentioning that problem $(P)$ is related to a class of interesting  applicative models. 
A common technique for the fabrication of semiconductor devices is the so-called {\it ion etching}, in which the material to be etched is bombarded with an ion beam (see \cite{F, R1,R2}). Mathematical modelling of the process gives rise to the Hamilton-Jacobi equation in one space dimension: 
$$
\left\{\begin{array}{ll}
U_t+ \varphi(U_x)=0 
& \quad\mbox{in}\  \R\times (0,T) \smallskip\\
U=U_0  &\quad\mbox{in}\  \R\times \{0\}\,,
\end{array}\right. \leqno{(HJ)}
$$ 
where $U=U(x,t)$ denotes the thickness  of the material and $\varphi$ is bounded, non-convex and vanishing at infinity. Clearly,  problem $(HJ)$ is related to $(P)$ by formally differentiating with respect to $x$ and setting $u=U_x$, $u_0=U_0'$. In this way space  discontinuous  solutions of $(HJ)$ correspond to Radon measure-valued solutions of $(P)$, which have a Dirac mass $\delta_{x_0}$ concentrated at any point $x_0$ where $U(\cdot,t)$ is discontinuous $(t\in(0,T))$.


\smallskip

Problem $(P)$ was studied in \cite{BSTT} under more general hypotheses on $\varphi$:
\begin{equation*} 
\varphi\in C([0,\infty))\,,\;\; \varphi(0)=0\,,\;\; \varphi'\in L^{\infty}(0,\infty)\,,\;\; \text{there exists 
\;$\displaystyle{ \lim_{u\to\infty}  \frac{\varphi(u)}{u}  =:C_{\varphi} \,.} $}
\leqno(A_1)
\end{equation*}
Without loss of generality one may assume that $C_\varphi=0$ (otherwise replace $x$ by $x-C_\varphi t$, see \cite{BSTT}).
If $u_0$ is any positive bounded Radon measure, an approximation approach can be used to construct
suitably defined entropy solutions of $(P)$ in a space of bounded Radon measures on $S$
(see Definitions \ref{deso}-\ref{enso} below and \cite[Theorem 3.2]{BSTT}; in the present section
we call such solutions ``constructed solutions''). However, an additional condition on solutions is needed for the well-posedness of $(P)$, since examples of nonuniqueness can be easily produced (see \cite{BSTT, DS}).
If $(H_0)$ and $(A_1)$ hold and $\varphi$ is bounded and {\it monotonic}, a uniqueness condition is known.
It prescribes the behaviour of the regular part  $u_r$ of the solution at points of the support of its singular part $u_s$:
\begin{equation}\label{crus}
(x_j,t)\in \supp u_s \ \Rightarrow \ 
\begin{cases}
{\rm ess} \lim_{x\to x_j^+}
u_r(x,t)=\infty &\text{if }\varphi'>0 \text{ in }[0,\infty) \\
{\rm ess} \lim_{x\to x_j^-} u_r(x,t)=\infty &\text{if }\varphi'<0\text{ in }[0,\infty)\,.
\end{cases}
\end{equation}
More precisely, in this case there exists at most one entropy solution of $(P)$, which satisfies \eqref{crus} and is strongly continuous at $t=0$ (\cite[Theorem 3.11]{BSTT}). If in addition
\begin{equation*}  \varphi\in C^3([0,\infty));\;
\exists\,L\ge-1, K\in \R\ \text{such that}\ \varphi''(u)[L\varphi(u)+K]\le -[\varphi'(u)]^2<0\,,
 \leqno(A_2)
\end{equation*}
every constructed entropy solution $u$ of problem $(P)$ satisfies \eqref{crus}, thus providing an existence and uniqueness theorem for $(P)$ (\cite[Theorem 3.12]{BSTT}). 

Observe that $(A_2)$ entails hat $\varphi$ is either increasing and concave or decreasing and  convex.
It is the aim of this paper to extend the above well-posedness results to the general case of a  bounded flux $\varphi$,  without assumptions about its monotonicity or convexity. To this purpose, we must find a general condition which replaces \eqref{crus}.


\subsection{A modified Cauchy problem }

Condition \eqref{crus} was suggested by the model problem (see \cite{BSTT})
\begin{equation}\label{ester}
\left\{\begin{array}{ll}
u_t+ \left[\varphi(u)\right]_x=0 
& \quad\mbox{in $S$}  \qquad \qquad \text{with $\varphi(u)=1-(1+u)^{-p}$ \quad ($p>0$)}
\smallskip\\
u=\delta_0  &\quad\mbox{in}\  \R\times \{0\}\,,
\end{array}\right. 
\end{equation}
where $T>1$. In fact, the unique constructed entropy solution of problem \eqref{ester} is 
\begin{equation}\label{sol1}
u_r(x,t):=\left[ (ptx^{-1})^{\frac{1}{1+p}}-1\right]\,\chi_A(x,t)\,,
\quad u_s(t):=\max\{1-t,0\}\delta_0 
\quad ((x,t)\in S\,),
\end{equation} 
where we have set
$$
A :=\{(x,t)\in S \, |\, 0<x\le p\,t, 0\le t < 1 \}\cup \{(x,t)\in S \, |\, \xi(t)\le x\le p\,t, 1\le t \le T  \},
$$ 
and $\xi(t)$ is defined by 
\begin{equation*}
\xi' = \dfrac{1-(pt\xi^{-1})^{-\frac{p}{1+p}}}{(pt\xi^{-1})^{\frac{1}{1+p}}-1} \quad\mbox{in $(1,T)$}, 
\qquad \xi(1)=0 \,.
\end{equation*}
Observe that the function  $u_r$ defined  in \eqref{sol1} diverges as $x\to 0^+$ if $t\in(0,1)$ - namely, as long as $u_s(\cdot,t)(\{0\})>0$,
in agreement with the first equality in \eqref{crus}. On the other hand, $u_r\equiv 0$ in the halfstrip $S_-:=(-\infty,0)\times (0,T)$, in particular $u_r(x,t)\to 0$ as $x\to 0^-$.

To generalize \eqref{crus} it is natural to address the ``modified Riemann problem'': 
\begin{equation}\label{rimo}
\left\{\begin{array}{ll}
u_t+ \left[\varphi(u)\right]_x=0 
& \quad\mbox{in $S
$} 
\smallskip\\
u=u_0  &\quad\mbox{in}\ \R\times \{0\}\,,
\end{array}\right. 
\end{equation}
with 
\begin{equation}\label{rimod}
u_{0r}=  u_- \chi_{(-\infty,0)} +u_+ \chi_{[0,\infty)}\,,\quad u_{0s}=\delta_0  \qquad (u_\pm \in[0,\infty))\,.\end{equation}

We seek entropy solutions of \eqref{rimo}-\eqref{rimod}. We introduce the sets $V_\pm \subseteq[0,\infty]$ of points which are {\em visible from the right}, 
\begin{equation}\label{deV}
\begin{aligned}
&V_+(\varphi):= \{u\in[0,\infty]\,|\, \varphi(u)\ge\varphi(s)\;\;\forall s\in(u,\infty] \}\\
& V_-(\varphi):= \{u\in[0,\infty]\,|\, \varphi(u)\le\varphi(s)\;\;\forall s\in(u,\infty] \} \,, 
\end{aligned}
\end{equation} 
and set 
\begin{equation*}
s_\pm (\varphi\,,u_0):=\inf \, V_\pm (\varphi)\cap  [u_0,\infty] \qquad \text{for }u_0\in [0,\infty)\,.
\end{equation*} 
For shortness we often write $V_\pm$ instead of $V_\pm(\varphi)$, and $s_{\pm}$ or $s_{\pm}(u_0)$ instead of $s_{\pm}(\varphi\,,u_0)$. Observe that for all $u_0\in [0,\infty)$
$$
\text{$s_+=\infty$,\, $s_-=u_0$ if  $\varphi$ is increasing}, \qquad \text{$s_+=u_0$,\, $s_-=\infty$ if  $\varphi$ is decreasing.}
$$

The following properties of $s_{\pm}$ are easily checked:
\begin{equation}\label{s-}
\begin{aligned}
& s_+=\inf \{u\in [u_0,\infty]\;|\,  \sup_{[u_0,u)} \varphi =\sup_{[u_0,\infty)}\varphi\}\,, \\
& s_-= \inf \{u\in [u_0,\infty]\;|\,  \inf_{[u_0,u)} \varphi =\inf_{[u_0,\infty)}\varphi\}\,;
\end{aligned} 
\end{equation}
\begin{equation}\label{svi}
\text{$s_\pm\in V_\pm$\,,  hence $s_\pm(\varphi\,,u_0)=\min \, V_\pm (\varphi)\cap  [u_0,\infty]$}\,.
\end{equation} 
Plainly, it follows from \eqref{s-}-\eqref{svi} that
\begin{equation}\label{ss}
s_{\pm}(s_{\pm}(u_0))=s_{\pm}(u_0)\,,
\end{equation}
\begin{equation}\label{mima}
\varphi(s_+(u_0))=  
\begin{cases}
\sup_{[s_+(u_0),\infty)}\varphi &\text{if $s_+(u_0)<\infty$} \\
\limsup_{s\to\infty} \varphi(s) &\text{if $s_+(u_0)=\infty$\,,} 
\end{cases} 
\end{equation}
\begin{equation}\label{mimabis}
 \varphi(s_-(u_0))=  
\begin{cases}
\inf_{[s_-(u_0),\infty)}\varphi\ &\text{if $s_-(u_0)<\infty$} \\
\liminf_{s\to\infty} \varphi(s) &\text{if $s_-(u_0)=\infty$\,,} 
\end{cases} 
\end{equation}
whence
\begin{equation}\label{p0-}
 \varphi'(s_{\pm}(u_0))=0\quad \text{if $u_0<s_{\pm}(u_0)<\infty$\,,} 
 \end{equation}
$$ 
\varphi' (s_+(u_0))\le0\quad \text{if $s_+(u_0)=u_0$\,,}\qquad 
\varphi' (s_-(u_0))\ge0\quad \text{if $s_-(u_0)=u_0$\,, }
$$
and
\begin{equation}\label{cfr}
\varphi(s_-(u_0))\le \liminf_{s\to\infty} \varphi(s) \le  \limsup_{s\to\infty} \varphi(s) \le \varphi(s_+(u_0))\quad \text{for $u_0\in[0,\infty]$\,.} 
 \end{equation}

Now consider the Riemann problems:
\begin{equation}\label{rv+}
\begin{cases}
v_t+[\varphi(v)]_x=0 &\text{in } S \\
v=   s_+(u_+) \chi_{(-\infty,0)} + u_+ \chi_{[0,\infty)}  &\text{in } \R\times\{0\}\,,
\end{cases}
\end{equation}
\smallskip
\begin{equation}\label{rv-}
\begin{cases}
v_t+[\varphi(v)]_x=0 &\text{in } S \\
v=   u_- \chi_{(-\infty,0)} +s_-(u_-) \chi_{[0,\infty)}  &\text{in } \R\times\{0\}
\end{cases}
\end{equation}
with $s_{\pm}(u_{\pm})<\infty$. 
Denote by $v_{\pm}$ the unique entropy solution of \eqref{rv+} and \eqref{rv-}, respectively. 
If $u_{\pm}=s_{\pm}(u_{\pm})$, there holds $v_{\pm}\equiv  s_{\pm}(u_{\pm})$ in $S$. 
On the other hand, if $u_{\pm}<s_{\pm}(u_{\pm})$, 
$v_{\pm}$ can be constructed in a standard way 
by considering the convex hull of $\varphi$ in the interval $[u_-,s_-(u_-)]$, respectively its concave hull 
in the interval $[u_+,s_+(u_+)]$, and the corresponding characteristics.
Plainly, by \eqref{p0-} there holds
\begin{equation}\label{vequis}
v_{\pm}= s_{\pm}(u_{\pm}) \quad\text{in } S_{\mp}\,,
\end{equation}
where $S_+:=(0,\infty)\times(0,T)$ and $S_-:=(-\infty,0)\times(0,T)$; moreover,
\begin{equation}\label{trari}
\exists \lim_{x\to 0^{\pm}}v_{\pm}(x,t)=: v_{\pm}(0^{\pm},t)\,, \quad\text{and} \quad v_{\pm}(0^{\pm},t)=s_{\pm}(u_{\pm}) \;\; \text{for all }t\in(0,T)\,.
\end{equation}

Now we can construct an entropy solution of problem \eqref{rimo}. Set
\begin{equation}\label{defini}
u_r(\cdot,t):=\begin{cases}
v_-(\cdot,t) &\text{in } (-\infty,0) 
\\
v_+(\cdot,t) &\text{in } (0,\infty) 
\end{cases} \,,\qquad
u_s(\cdot,t):=A(t)\delta_0\,, 
\end{equation}
where 
\begin{equation}\label{A(t)}
A(t):=1-[\varphi(s_+(u_+))-\varphi(s_-(u_-))]\,t\,. 
\end{equation}
Observe that $\varphi(s_+(u_+))\ge\varphi(s_-(u_-))$ by \eqref{cfr}. If $\varphi(s_+(u_+))>\varphi(s_-(u_-))$,  the measure  $u$ defined by \eqref{defini}-\eqref{A(t)} is positive on $\R\times(0,\tau)$  with
$$
\tau:=\frac 1{\varphi(s_+(u_+))-\varphi(s_-(u_-))}\,.
$$
It is easily seen that $u$ is an entropy solution of the modified Riemann problem \eqref{rimo} in 
$\R\times(0,\min\{\tau,T\})$ (see Definitions \ref{deso}-\ref{enso}). 
By \eqref{trari}-\eqref{defini},  there holds
\begin{equation}\label{traribis}
u_r(0^{\pm}\!,t):=\lim_{x\to 0^{\pm}}u_r(x,t)=s_{\pm}(u_{\pm}) \quad\text{for all $t\in(0,\min\{\tau,T\})$}\, .
\end{equation}

If $\tau\ge T$, the result follows. Otherwise, we set $ u_s(\cdot,t):=0$ for all $t\in(\tau,T]$ and  continue the solution  in $(\tau,T]$, with initial data $u_r(\cdot,\tau)$, using the standard theory of scalar conservation laws. 
If $\varphi(s_+(u_+))=\varphi(s_-(u_-))$, it is easily seen that $u$ is an equilibrium solution in $S$:
$A(t)\equiv1$ in $[0,T]$, and, by \eqref{cfr} and the definition of $s_\pm$, 
$\varphi\equiv \varphi(s_+(u_+)=\varphi(s_-(u_-)$ in the interval $(\min\{u_-,u_+\},\infty)$, thus $u_r$ is constant in $S$ (see \eqref{defini}). One easily generalizes the above discussion to the case that $s_{\pm}(u_{\pm})=\infty$.

\smallskip

It is worth revisiting problem \eqref{ester} in the light of the above remarks. Since in this case $u_{\pm}=0$ and $\varphi$ is increasing, there holds $s_+(u_+)=\infty$, $s_-(u_-)=0$, $\varphi(s_+(u_+))=1$ and $\varphi(s_-(u_-))=0$, whence 
(see \eqref{traribis}) 
\begin{equation*}
\lim_{x\to 0^+}u_r(x,t)=\infty\,, \quad \lim_{x\to 0^-}u_r(x,t)=0\,, \quad \text{$A(t)=1-t$\,\quad for $t\in[0,1]$\,,}
\end{equation*}
in agreement with \eqref{sol1}.


\subsection{Compatibility conditions } To address problem $(P)$ under assumption $(H_0)$ we need a more general condition
than \eqref{traribis}, which is only suitable for the modified Riemann problem. To this purpose, observe that equalities \eqref{ss} and \eqref{traribis} entail
\begin{equation}\label{trariter}
u_r(0^{\pm}\!,t) =s_{\pm}(u_r(0^{\pm}\!,t))\,
\end{equation}
for all $t\in(0,T)$ such that $u_s(\cdot,t)>0$. For problem \eqref{ester}  
the equality at $0^+$ coincides with the first equality in \eqref{crus}, while that at $0^-$ is trivially satisfied. So, if $\{0\}\times(0,t)\in \supp u_s$, it is natural to regard  \eqref{trariter} as the desired generalization of \eqref{crus}.

\smallskip

Set $H_-(u):=- \chi_{(-\infty,0)}(u)$ $(u\in\R)$. It is easily seen that condition \eqref{trariter} can be rephrased as 
\begin{equation}\label{trariqua}
\begin{aligned}
&  H_-(u_r(0^+\!,t)-k)[\varphi(u_r(0^+\!,t))-\varphi(k)]\le0
\\
& H_-(u_r(0^-\!,t)-k)[\varphi(u_r(0^-\!,t))-\varphi(k)]\ge0
\end{aligned}
\qquad\text{for all $k\in[0,\infty)$\,.}
\end{equation}
Formally \eqref{trariqua} is equivalent to the {\em compatibility condition}  
\begin{equation}\label{trariqui}
\begin{aligned}
&  [\sgn(u_r(0^+,t)-k) - \sgn(a_0(t)-k)][\varphi(u_r(0^+\!,t))-\varphi(k)]\le0\,
\\
& [\sgn(u_r(0^-,t)-k) - \sgn(a_0(t)-k)][\varphi(u_r(0^-\!,t))-\varphi(k)]\ge0
\end{aligned}
\end{equation}
between the traces $u_r(0^{\pm}\!,t)$ and the boundary data $a_0(t)=\infty$, for all $k,t$ as above. 
It was shown in \cite{BLN, Te} that the initial-boundary value problems 
\begin{equation*}
\begin{cases}
v_t+[\varphi(v)]_x=0 &\text{in } S_+ \\
v =   a_0  &\text{in } \{0\}\times (0,T)\\
v=    u_+   &\text{in } [0,\infty)\times\{0\}\,,
\end{cases}
\qquad 
\begin{cases}
v_t+[\varphi(v)]_x=0 &\text{in } S_- \\
v =   a_0  &\text{in } \{0\}\times (0,T)\\
v=    u_-   &\text{in } (-\infty,0]\times\{0\}\,
\end{cases}
\end{equation*}
are well posed, if $a_0\in BV(0,T)$ and \eqref{trariqui} holds. 
This gives an alternative interpretation of the construction used to solve the modified Riemann problem \eqref{rimo} (see  \eqref{rv+}-\eqref{rv-}): as long as the Dirac delta at $x=0$ survives, it behaves like a barrier which decouples the evolution of the regular part of the solution on either side of the singularity, imposing the two Dirichlet conditions $u_r(0^\pm,t)=\infty$ at $x=0$. The evolution of the delta at $t=0$ is completely determined by local mass exchange through $x=0$. 

\smallskip

The above considerations suggest a constructive approach to address problem $(P)$ under assumption $(H_0)$. By the results in \cite{BSTT} there is a positive time $\tau$ until which all singularities persist, thus the real line is the disjoint union of $p+1$ intervals. In each interval we solve the initial-boundary value problem for the conservation law in $(P)$, the initial data being the restriction of $u_{0r}$ to that interval, with ``boundary conditions equal to infinity" - or, equivalently, by imposing the analogue of \eqref{trariqua} to be satisfied at each point $x_j$, $j=1,\dots,p$. The function determined  by this procedure is, by definition, the regular part of a Radon measure, whose singular part is defined in analogy with \eqref{defini}-\eqref{A(t)}. It is proven that this measure is the unique entropy solution of $(P)$ (in the sense of Definitions \ref{deso}-\ref{enso}) until the time $t=\tau$. If $\tau<T$ we iterate the procedure in $\R\times(\tau,T)$ with a smaller number a singularities, thus well-posedness of $(P)$ follows in a finite number of steps (see Theorem \ref{exiuni}).

A technical obstruction to the above program is that the solution constructed in each interval need not have traces at the points $x_j$. This difficulty is overcome by using a weak analogue of condition \eqref{trariqua} (see \eqref{coco}) and the $L^\infty$-theory of initial-boundary value problems developed in \cite{O}.

By the finite speed of propagation of solutions of hyperbolic conservation laws, uniqueness proofs
are local in space. Then it can be easily checked that our results remain valid, if condition $(H_0)$ is relaxed
 to the case that $u_{0s}$ is a {\it locally} finite superposition of Dirac masses (namely, in every bounded interval the number of Dirac masses is finite). The case of more general $u_{0s}$ is open.

The paper is organized as follows.  After recalling some preliminaries (see Section~\ref{preli}), the main results of the paper are presented in Section~\ref{resu}, whereas  Sections~\ref{prexi}-\ref{comp} are devoted to their proofs.


\section{Preliminaries}\label{preli}
\setcounter{equation}{0}

Let $\chi_E$ denote the characteristic function of $E\subseteq\R$. For all $u\in\R$, we set $u_\pm= \max \{\pm u,0\}$, 
$H_\pm(u)=\pm \chi_{\{\pm u>0\}}(u)$, $\sgn (u)=H_+(u)+H_-(u)$. For every real function $f$ on $\R$ and $x_0\in\R$ we say that 
$$
\text{{\rm ess\,}$\lim_{x\to x_0^{\pm}}  f(x) =l\in\R$\,,}
$$ 
if there is a  null set $E^*\subseteq\R$  such that $f(x_n)\to l$ if $\left\{x_n\right\}\subseteq\R\setminus\! (E^*\!\cup\!\{x_0\})$, $x_n\to x_0^{\pm}$. 

For every open subset $\Omega\subseteq\R$ we denote by $\mathcal{M}(\Omega)$ the space of Radon measures on $\Omega$,  by $\mathcal{M}^+(\Omega)$  the cone of its nonnegative elements. If $\mu,\nu\in\mathcal{M}(\Omega)$, we say that $\mu\le\nu$ in $\mathcal{M}(\Omega)$ if $\nu-\mu\in\mathcal{M}^+(\Omega)$. We denote by $C_c(\Omega)$ the space  of continuous real functions with compact support in $\Omega$, and by $\lla \cdot, \cdot\rra_{\Omega}$ the duality map between $\mathcal M(\Omega)$  and  $C_c(\Omega)$. A sequence $\{\mu_n\}$ of Radon measures on $\R$ converges weakly* to a Radon measure $\mu$, $\mu_n\stackrel{*}\rightharpoonup \mu$, if $\left\langle \mu_n,\rho\right\rangle_{\R}\to \left\langle \mu,\rho\right\rangle_{\R}$ for all $\rho \in C_c(\R)$. For any compact $K\in\R$ the space $\mathcal{M}(K)$ is a Banach space with norm $\|\mu \|_{\mathcal M(K)}:=|\mu|(K)$.
A sequence $\{\mu_n\}$ converges strongly to $\mu$ in $\mathcal{M}(K)$ if  $\|\mu_n-\mu\|_{\mathcal{M}(K)}\to 0$ as  $n\to\infty$.
Similar definitions are used for Radon measures on any subset of $S:=\R\times (0,T)$. 

Every $\mu\in\mathcal{M}(\R)$ has a 
unique decomposition $\mu=\mu_{ac}+\mu_s$, with $\mu_{ac} \in\mathcal{M}(\R)$
absolutely continuous and $\mu_s\in\mathcal{M}(\R)$ singular with respect to the Lebesgue  measure. 
We  denote by $\mu_r\in L^1_{loc}(\R)$ the density of $\mu_{ac}$. Every function $f\in L^1_{loc}(\R)$ can be identified to an absolutely continuous Radon measure on $\R$; we shall denote this measure by the same symbol $f$ used for the function. 

The restriction $\mu \, \lefthalfcup E$ of $\mu\in \mathcal{M}(\R)$ to a Borel set $E\subseteq \R$ is defined by 
$(\mu \, \lefthalfcup E) (A):=\mu(E\cap A)$ for any Borel set $A\subseteq\R$. 
Similar notations are used for $\mathcal{M}(S)$. 

\smallskip

 We shall use measures $u\in \mathcal{M}(S)$ which, roughly speaking, admit a parametrization with respect to the time variable:
\begin{definition}\label{dli}
We denote by $L^{\infty}(0,T;\mathcal{M}^+(\R))$ the set of nonnegative Radon measures $u\in \mathcal{M}^+(S)$ such that 
for a.e.~$t\in (0,T)$ there is a measure $u(\cdot,t)\in \mathcal{M}^+(\R)$ with the following properties:

\noindent $(i)$ if $\zeta\in C([0,T];C_c(\R))$ the map $t\mapsto \left\langle u(\cdot,t),\zeta(\cdot,t)\right\rangle_{\R}$ 
belongs to $L^1(0,T)$ and
\begin{equation}\label{eq.disintegrazioneU}
\left\langle u,\zeta\right\rangle_{S}=\int_0^T\left\langle u(\cdot,t),\zeta(\cdot,t)\right\rangle_{\R}\,dt\,;
\end{equation}

\noindent $(ii)$ the map $t\mapsto \|u(\cdot,t)\|_{\mathcal{M}(K)}$ belongs to $L^{\infty}(0,T)$ for every compact $K\in\R$.
\end{definition} 

\begin{remark} Definition \ref{dli} implies that for all  $\rho\in C_c(\R)$ the map  $t\mapsto \left\langle u(\cdot,t),\rho\right\rangle_{\R}$ 
is measurable, thus the map   $u:(0,T) \to \mathcal{M}(\R)$ is weakly* measurable. For simplicity we prefer the notation $L^{\infty}(0,T;\mathcal{M}(\R))$ to the more correct one $L^{\infty}_{w*}(0,T;\mathcal{M}(\R))$. \end{remark}

Observe that $u_r\in L^{\infty}(0,T;L^1_{loc}(\R))$ if $u\in L^{\infty}(0,T;\mathcal{M}^+(\R))$. Conversely, every nonnegative 
$f\in L^{\infty}(0,T;L^1_{loc}(\R))$ defines a measure belonging to 
$ L^{\infty}(0,T;\mathcal{M}^+(\R))$.

By $C([0,T];\mathcal{M}^+(\R))$ we denote the subset of strongly continuous mappings from $[0,T]$ into $\mathcal{M}^+(\R)$ - namely, $u\in C([0,T];\mathcal{M}^+(\R))$ if for all $t_0\in[0,T]$ and for every compact $K\in\R$ there holds $ \|u(\cdot,t)-u(\cdot,t_0)\|_{\mathcal{M}(K)}\to0$ as $t\to t_0$.

If $u\in L^{\infty}(0,T;\mathcal{M}^+(\R))$, also $u_{ac}, u_s\in L^{\infty}(0,T;\mathcal{M}^+(\R))$ and, by  \eqref{eq.disintegrazioneU},  
\begin{equation}\label{disicomp}
\lla u_{ac}\,, \zeta\rra_{S}\,=\iint_{S} u_r \,\zeta\,dxdt \,,
\quad 
\lla u_s, \zeta \rra_{S} \,=\int_0^T\!\! \lla u_s(\cdot,t),\zeta(\cdot,t) \rra_{\R} dt
\end{equation}
if $\zeta\in C([0,T];C_c(\R))$. One easily checks that for a.e.~$t\in (0,T)$
\begin{equation}\label{us(t)=u(t)s}
u_{ac}(\cdot,t)=[u(\cdot,t)]_{ac}\,, \quad u_s(\cdot,t)=[u(\cdot,t)]_s\,, \quad u_r(\cdot,t)=[u(\cdot,t)]_r\,,
\end{equation}
where $[u(\cdot,t)]_r$ denotes the density of the measure $[u(\cdot,t)]_{ac}$: for   $\rho\in C_c(\R)$ 
$$
\lla [u(\cdot,t)]_{ac}, \rho\rra_{\R}  =\int_{\R} [u(\cdot,t)]_r \,\rho \,dx 
=\int_{\R} u_r(\cdot,t)\, \rho \,dx \quad\text{for a.e.~$t\in(0,T)$.} 
$$  
In view of  \eqref{disicomp}-\eqref{us(t)=u(t)s}, we shall always  identify the quantities which appear 
on either side of equalities \eqref{us(t)=u(t)s}.


\section{Results}\label{resu}
\setcounter{equation}{0}

For any open $\Omega\subseteq\R$ and $\tau\in(0,T)$ set $Q_{\tau}:=\Omega\times (0,\tau)$. Solutions of problem $(P)$ are meant in the following sense. 
\begin{definition}\label{deso}  
A measure $u\in L^\infty(0,T;\mathcal{M}^+(\Omega))$ is called a {\it solution} of problem $(P)$ in $Q_{\tau}$  if for all $\zeta\in C^1([0,\tau];C^1_c(\Omega))$, $\zeta(\cdot,\tau)=0$ in $\Omega$ there holds
\begin{equation}\label{ewf}
\iint_{Q_{\tau}} \big[u_r\zeta_t+ \varphi(u_r) \,\zeta_x\big]\,dxdt+\int_0^{\tau}\lla u_s(\cdot,t), 
\zeta_t(\cdot,t)\rra_{\Omega}dt= - \lla u_0,\zeta(\cdot,0)\rra_{\Omega}\,.
\end{equation}
Solutions of $(P)$ in $S$ are simply referred to as ``solutions of $(P)$''.
\end{definition}

\begin{definition}\label{enso}
A solution  of $(P)$ in $Q_{\tau}$ is called an {\em entropy solution in $Q_{\tau}$} if it satisfies the {\em entropy inequality}
\begin{eqnarray}\label{mkru}
&&\iint_{Q_{\tau}} \left\{|u_r-k|\,\zeta_t+\sgn(u_r-k)\left [\varphi(u_r)-\varphi(k)\right ]\zeta_x\right\}dxdt 
+\\
&+& \int_0^\tau\left\langle u_s(\cdot,t),\zeta_t(\cdot,t)\right\rangle_{\Omega}\,dt \geq - 
\int_{\Omega} |u_{0r}(x)-k|\,\zeta(x,0)
\,dx - \left\langle u_{0s}, \zeta(\cdot,0)\right\rangle_{\Omega} \nonumber  
\end{eqnarray}
 for all $\zeta\in C^1([0,\tau];C^1_c(\Omega))$, $\zeta\geq0$, $\zeta(\cdot,\tau)=0$ in $\Omega$, and for all $k\in[0,\infty)$. 
\end{definition}
\begin{remark}\label{susu}  
{\em Entropy subsolutions and supersolutions}  of $(P)$ in $Q_{\tau}$ are defined by requiring the following inequalities to be satisfied: 
\begin{eqnarray}\label{subkru}
&&\iint_{Q_{\tau}} \left\{[u_r-k]_+\,\zeta_t+H_+(u_r-k)\left [\varphi(u_r)-\varphi(k)\right ]\zeta_x\right\}dxdt 
+\\
&+&\int_0^\tau\left\langle u_s(\cdot,t),\zeta_t(\cdot,t)\right\rangle_{\Omega
}\,dt \geq - \int_{\Omega
} [u_{0r}-k]_+\,\zeta(x,0)\,dx - \left\langle u_{0s}, \zeta(\cdot,0)\right\rangle_{\Omega
}, \nonumber  
\end{eqnarray}
respectively
\begin{equation}\label{superkru}
\iint_{Q_{\tau}} \!\!\left\{[u_r\!-\!k]_-\,\zeta_t\!+\!H_-(u_r\!-\!k)\left [\varphi(u_r)\!-\!\varphi(k)\right ]\zeta_x\right\}dxdt 
\ge - \!\int_{\Omega} [u_{0r}\!-\!k]_-\,\zeta(x,0)\,dx 
\end{equation}
 for all $\zeta$ and $k$ as above. It is easily seen that $u$ is an entropy solution if and only if it is both an entropy subsolution and an entropy supersolution. 
  \end{remark}

Let $\beta\in C^1_c(0,T)$, $\beta\ge 0$, and $k\in[0,\infty)$. We shall prove below (see Lemma \ref{lim}) that, if $(H_0)$-$(H_1)$ hold, for every entropy solution of $(P)$ the limits
\begin{equation}\label{defili}
{\rm ess} \lim_{x\to x_j^{\pm}} \int_0^T H_-(u_r(x,t)-k)[\varphi(u_r(x,t))-\varphi(k)]\beta(t)\,dt \qquad (j=1,\dots,p
) 
\end{equation}
exist and are finite. It is also known that, if $(H_0)$-$(H_1)$ are satisfied, $j=1,\dots,p$ and $u$ is a solution of problem $(P)$,
 \begin{equation}\label{dewait}
 \text{$\forall x_j$  $\exists  t_j\in(0,T]$ such that }
\begin{cases}
u_s(\cdot,t)(\{x_j\})>0 &\text{for a.e. $t\in[0,t_j)$} \\
u_s(\cdot,t)(\{x_j\})=0 &\text{for a.e. $t\in(t_j,T)$} 
\end{cases} 
\end{equation}
 (see \cite[Theorem 3.5]{BSTT}). Then we can state the following definition.
\begin{definition}\label{decoco}
Let $(H_0)$-$(H_1)$ be satisfied, let $j=1,\dots,p$ and let $\tau\in(0,t_j]$.
An entropy solution of $(P)$ satisfies the {\em compatibility condition at $x_j$ in} $[0,\tau]$ 
if for all $\beta\in C^1_c(0,\tau)$, $\beta\ge 0$, and $k\in[0,\infty)$
\begin{equation}\label{coco}
\begin{aligned}
& {\rm ess} \lim_{x\to x_j^+} \int_0^{\tau} H_-(u_r(x,t)-k)[\varphi(u_r(x,t))-\varphi(k)]\beta(t)\,dt \le 0\,, \\
& {\rm ess} \lim_{x\to x_j^-} \int_0^{\tau} H_-(u_r(x,t)-k)[\varphi(u_r(x,t))-\varphi(k)]\beta(t)\,dt \ge0\,.
\end{aligned}
\end{equation}
\end{definition}
Now our main result can be stated as follows.
\begin{theorem}\label{exiuni}
 Let $(H_0)$-$(H_1)$ be satisfied. Then there exists  a unique entropy solution of problem $(P)$ which belongs to $C([0,T];\mathcal{M}^+(\R))$ and satisfies the compatibility condition at $x_j$ in $[0,t_j]$ for all  $j=1,\dots,p$.
\end{theorem}
According to Theorem \ref{exiuni}, the compatibility condition defines a well-posedness class for entropy solutions of $(P)$ under assumptions $(H_0)$-$(H_1)$.

\smallskip

We shall also prove a comparison result for solutions of $(P)$ whose initial data satisfy assumption $(H_0)$: 
 \begin {theorem}\label{compa}  Let $(H_1)$ be satisfied. Let $v_0\in\mathcal{M}^+(\R)$ satisfy $(H_0)$, and let $u_0\le v_0$ in $\mathcal{M}(\R)$. Let $u,v\in C([0,T];\mathcal{M}^+(\R))$ be the unique entropy solutions of $(P)$ with initial data $u_0, v_0$ given by Theorem \ref{exiuni}. Then there holds $u(\cdot,t)\le v(\cdot,t)$  in $\mathcal{M}(\R)$ for all $t\in[0,T]$. 
 \end{theorem}


\section{Proof of existence}\label{prexi} 
\setcounter{equation}{0}

In this section we prove the existence part of Theorem \ref{exiuni}: 
\begin{theorem}\label{exi}
 Let $(H_0)$-$(H_1)$ be satisfied. Then there exists  an entropy solution of problem $(P)$ which satisfies the compatibility condition at $x_j$ in $[0,t_j]$ for all  $j=1,\dots,p$. Moreover, $u$ belongs to $C([0,T];\mathcal{M}^+(\R))$.
\end{theorem}
To prove Theorem \ref{exi} we need some preliminary results.


\subsection{Preliminary results}

\begin{lem}\label{lim}
Let $u$ be an entropy supersolution of $(P)$, and let $\beta\in C^1_c(0,T)
$, $\beta\ge 0$. Then:

\noindent $(i)$ for every $k\in[0,\infty)$ the  distributional derivative of the function 
\begin{equation}\label{hacca}
x\mapsto  -\int_0^T\!\! H_-(u_r(x,t)-k)\left [\varphi(u_r(x,t))-\varphi(k)\right ]\beta(t)\,dt + kT \|\beta'\|_\infty\, x
\end{equation}
is nonnegative; 
 
 \noindent $(ii)$  for every $x_0\in\R$ and $k\in[0,\infty)$  the limits
 \begin{equation}\label{ali}
{\rm ess} \lim_{x\to x_0^{\pm}} \int_0^T H_-(u_r(x,t)-k)[\varphi(u_r(x,t))-\varphi(k)]\beta(t)\,dt 
\end{equation}
exist and are finite.
\end{lem}

\begin{proof}
Let $\alpha\in C^1_c(\R)$, $\alpha\ge0$. Choosing $\zeta(x,t)=\alpha(x)\beta(t)$ in \eqref{superkru} with $Q_{\tau}=S$ gives
$$
\iint_S \! \! \left\{[u_r(x,t)\! -\! k]_-\,\alpha(x)\beta'(t)\! +\! H_-(u_r(x,t)\! -\! k)\! \left [\varphi(u_r(x,t))\! -\! \varphi(k)\right ]\!  \alpha'(x)\beta(t)\right\}\! dxdt\!  \ge \! 0.
$$
Since $0\le [u_r-k]_-\le k$, from the above inequality we get 
$$
-\! \int_{\R} \! \left(\int_0^T \! \! H_-(u_r(x,t)\! -\! k)\! \left [\varphi(u_r(x,t))\! -\! \varphi(k)\right ]\! \beta(t)\,dt\right) \, \alpha'(x)dx \le kT  \|\beta'\|_\infty\! \!  \int_{\R}\!  \alpha(x)dx,
$$
whence claim $(i)$ follows. 

Therefore, the  distributional derivative of function \eqref{hacca} is a Radon measure. Clearly, the same holds  for the  distributional derivative, say $\mu$, of the function $\mathcal{H}\in L^1_{loc}(\R)$, 
\begin{equation*}
\mathcal{H}(x):= -\int_0^T\!\! H_-(u_r(x,t)-k)\left [\varphi(u_r(x,t))-\varphi(k)\right ]\beta(t)\,dt \,.
\end{equation*}
Fix any $\bar x\in\R$ and set $f_\mu(x):=\mu((\bar x, x])$ if $x\ge \bar x$, $f_\mu(x):=-\mu((x,\bar x])$ if $x< \bar x$. Then $f_\mu$ is continuous from the right, and coincides a.e. with $\mathcal{H}$ on every compact $K\subset\R$ up to a constant, possibly depending on $K$ ($e.g.$, see \cite[Theorem 3.28]{AFP}). Hence the claim follows.
\end{proof}

In the following we set $I_1:=(-\infty,x_1)$, $I_j:=(x_{j-1},x_j)$ for $j=2,\dots,p$, $I_{p+1}:=(x_p,\infty)$ and $S_j:=I_j\times (0,T)$ for $j=1,\dots,p+1$.
\begin{lem}\label{lime} Let   $(H_0)$-$(H_1)$ hold, and let $u$ be an entropy solution of $(P)$. Then for all $\beta\in C^1_c(0,T)$, $\beta\ge 0$, $k\in[0,\infty)$ and $j=1,\dots,p$ the limits 
\begin{equation}\label{alibis}
{\rm ess} \lim_{x\to x_j^{\pm}} \int_0^T \sgn(u_r(x,t)-k)[\varphi(u_r(x,t))-\varphi(k)]\beta(t)\,dt
\end{equation}
exist and are finite.
\end{lem}
\proof
 We only prove the claim for the limit from the right,  the proof being similar for the other. Let $j=1,\dots,p$ be fixed. Since $\varphi$ is bounded, by \cite[Proposition 3.3]{BSTT} the singular part  of every entropy solution of $(P)$ is nonincreasing in time, hence $u_s(\cdot,t)(I_{j+1})=0$ for any $t\in [0,T]$. 
 Let $\alpha\in C^1_c(I_{j+1})$, $\alpha\ge0$. 
Choosing  $\zeta(x,t)=\alpha(x)\beta(t)$ in \eqref{mkru} with $Q_{\tau}=S$ gives
$$ 
\iint_{S_{j+1}}\!\!\!
 \{|u_r(x,t)\!-\!k|\,\alpha(x)\beta'(t)\!+\sgn(u_r(x,t)\!-\!k)\!\left [\varphi(u_r(x,t))\!-\!\varphi(k)\right ]\alpha'(x)\beta(t)\}dxdt\! \ge\! 0\,.
$$ 
Since $0\le |u_r-k|\le u_r+ k$, we have
\begin{eqnarray*}
&&-\int_{I_{j+1}}
\!\!dx\,\alpha'(x) \left(\int_0^T \!\!\!\sgn(u_r(x,t)-k)\left [\varphi(u_r(x,t))-\varphi(k)\right ]\beta(t)\,dt\right)\,  \le \\
&\le& \!\! \|\beta'\|_\infty \int_{I_{j+1}}
\!\!dx\,\alpha(x)\left (\int_0^Tu_r(x,t)\,dt + kT \right)\,= \\
&=&\!\! -\; \|\beta'\|_\infty \int_{I_{j+1}}
\!\!dx\,\alpha'(x)\left (\int_0^T\!\!\int_{x_j}^x u_r(y,t)\,dydt + kTx \right)\,.
\end{eqnarray*}
The above inequality implies that the  distributional derivative
 of the map 
\begin{eqnarray*}
x\mapsto \!\!&-&\!\!\!\!\int_0^T \!\!\sgn(u_r(x,t)-k)\left [\varphi(u_r(x,t))-\varphi(k)\right ]\beta(t)\,dt \,+ \\
&+&\!\!\! \|\beta'\|_\infty\left (\int_0^T\!\!\int_{x_j}^x u_r(y,t)\,dydt + kTx \right)
\end{eqnarray*}
is nonnegative in $I_{j+1}$. Arguing as in the proof of Lemma \ref{lim} the claim follows.
\qed
\begin{lem}\label{incom}
Let   $(H_0)$-$(H_1)$ hold, and let $u$ be an entropy solution of $(P)$. Then for every $j=1,\dots,p$:

\smallskip

\noindent $(i)$  there exist $h_j^-, h_j^+\in L^{\infty}(0,T)$, $h_j^{\pm}\ge0$ such that for all $\beta\in C^1_c(0,T)$ 
\begin{equation}\label{lidesi}
{\rm ess} \lim_{x\to x_j^{\pm}}\int_0^T \varphi(u_r(x,t))\beta(t)\,dt=\int_0^T h_j^{\pm}(t)\beta(t)\,dt\,; 
\end{equation}

\noindent $(ii)$ if $u$ satisfies the compatibility condition \eqref{coco} at $x_j$ in $[0,\tau]$,  
 there holds
\begin{equation}\label{ext_h2}
h_j^-\le \liminf_{k\to\infty}\,\varphi(k)\le \limsup_{k\to\infty}\,\varphi(k) \le h_j^+ \quad \text{a.e.~in }(0,\tau)\,. 
\end{equation}
\end{lem}
\begin{remark}
By standard density arguments, from \eqref{lidesi} we get
\begin{equation}\label{lidesi1}
{\rm ess} \lim_{x\to x_j^{\pm}}\int_0^T \varphi(u_r(x,t))\zeta(x,t)\,dt=\int_0^T h_j^{\pm}(t)\zeta(x_j,t)\,dt 
\end{equation}
for every $\zeta\in L^1(0,T;C_c(U_j))$ with $x_j\in U_j\subseteq\R$, $U_j$ open.
\end{remark}

\noindent {\em Proof of Lemma \ref{incom}.} $(i)$ We only prove the limit from the right. Since $\sgn u=1+2H_-(u)$ for $u\in\R$, by \eqref{ali}-\eqref{alibis} the limit in the left-hand side of \eqref{lidesi} exists and is finite. On the other hand, for every sequence $\{x_n\}$ converging to $x_j^+$ the sequence $\{\varphi(x_n)\}$ is bounded in $L^{\infty}(0,T)$, hence there exist a subsequence $\{x_{n_k}\}\subseteq \{x_n\}$ and a function $h_j^+\in L^{\infty}(0,T)$ such that $\varphi(x_{n_k})\stackrel{*}\rightharpoonup h_j^+$ in $L^{\infty}(0,T)$. 

\smallskip

\noindent $(ii)$ We only prove the last inequality in \eqref{ext_h2}. Since  $u$ is a  solution of $(P)$ in $I_{j+1}\times(0,\tau)$, by \eqref{ewf} there holds
\begin{equation*}
\int_0^{\tau}\!\!\!\int_{I_{j+1}}
 \big\{ (u_r-k)\xi_t+ [\varphi(u_r)-\varphi(k)] \,\xi_x\big\}\,dxdt= - \int_{I_{j+1}}[ u_{0r}(x)-k]\,\xi(x,0)\,dx
\end{equation*}
for all $k\in[0,\infty)$  and $\xi\in C^1([0,\tau];C^1_c(I_{j+1}))$, $\xi(\cdot,\tau)=0$ in $I_{j+1}$. Let 
\begin{equation}\label{defeta}
\eta_{\epsilon}(x):=\frac{2(x-x_j)-\epsilon}{\epsilon}\chi_{[x_j+\epsilon/2,x_j+\epsilon]}(x)+\chi_{(x_j+\epsilon,x_{j+1}]}(x)\qquad (x\in I_{j+1})
\end{equation}
and let $\zeta\in C^1([0,\tau];C^1_c([x_j,x_{j+1})))$, $\zeta(\cdot,\tau)=0$ in $I_{j+1}$ (here $x_{j+1}=\infty$ if $j=p$). By standard arguments we can choose 
$\xi=\zeta\eta_{\epsilon}$
in the above equality, and obtain
$$
\begin{aligned}
&\int_0^{\tau}\!\!\!\!\int_{I_{j+1}}  \!\!\!\!\big\{ (u_r-k)\zeta_t\eta_{\epsilon}+ [\varphi(u_r)-\varphi(k)] \,\zeta_x\eta_{\epsilon}\big\}dxdt
+\!\! \int_{I_{j+1}}\!\!\![ u_{0r}(x)-k]\,\zeta(x,0)\eta_{\epsilon}(x)dx=\\
=&\; -\frac2\epsilon\int_0^{\tau}\!\!\!\int_{x_j+\epsilon/2}^{x_j+\epsilon} [\varphi(u_r)-\varphi(k)]\,\zeta\,dxdt\,.
\end{aligned}
$$
Letting $\epsilon\to 0^+$ in the above equality plainly gives (see \eqref{lidesi1}): 
\begin{eqnarray}\label{mkjn1}
&&\int_0^{\tau}\!\!\!\int_{I_{j+1}}  \big\{ (u_r-k)\zeta_t+ [\varphi(u_r)-\varphi(k)] \,\zeta_x\big\}\,dxdt+ \int_{I_{j+1}}[ u_{0r}(x)-k]\,\zeta(x,0)\,dx\,= \\
&=& -\;{\rm ess}\lim_{x\to x_j^+}\int_0^{\tau}\![\varphi(u_r(x,t))-\varphi(k)] \,\zeta(x,t)\,dt
= -\int_0^{\tau}\![h_j^+(t)-\varphi(k)] \,\zeta(x_j,t)\,dt\,.  \nonumber
\end{eqnarray}

Since  $u$ is an entropy solution of $(P)$ in $I_{j+1}\times(0,\tau)$, arguing as before we obtain
$$
\begin{aligned}
&\int_0^{\tau}\!\!\!\!\int_{I_{j+1}}\!\!\!\!\!\left\{|u_r-k|\,\zeta_t+\sgn(u_r-k)\left [\varphi(u_r)-\varphi(k)\right ]\zeta_x\right\}\!dxdt 
+ \!\!\!\int_{I_{j+1}}\!\!\!|u_{0r}(x)-k|\,\zeta(x,0)dx\! \ge \\
\ge&\; -\,{\rm ess}\lim_{x\to x_j^+}\int_0^{\tau}\sgn(u_r(x,t)-k)\left[\varphi(u_r(x,t))-\varphi(k)\right ]\zeta(x,t)\,dt
\end{aligned}
$$
for all  $\zeta$ as above,  $\zeta\ge0$. Choosing $\zeta(x,t)=\alpha(x) \beta(t)$ with  $\alpha\in C^1_c([x_j,x_{j+1}))$, $\alpha\ge0$ and $\beta\in C^1([0,\tau])$, $\beta\ge0$, $\beta(\tau)=0$, by the compatibility condition \eqref{coco} there holds:
\begin{eqnarray}\label{mkjn3}
&&\;\;\int_0^{\tau}\!\!\!\int_{I_{j+1}} \left\{|u_r-k|\,\zeta_t+\sgn(u_r-k)\left [\varphi(u_r)-\varphi(k)\right ]\zeta_x\right\}dxdt \,+ \\
&&+ \int_{I_{j+1}}|u_{0r}(x)-k|\,\zeta(x,0)\,dx\, + {\rm ess}\lim_{x\to x_j^+}\int_0^{\tau}\! \left[\varphi(u_r(x,t))-\varphi(k)\right ]\zeta(x,t)\,dt\,\ge  \nonumber  \\
&&\ge -2\,{\rm ess}\lim_{x\to x_j^+}\int_0^{\tau}H_-(u_r(x,t)-k)\left[\varphi(u_r(x,t))-\varphi(k)\right ]\zeta(x,t)\,dt\,=\nonumber \\
&&=-2\,\alpha(x_j)\,{\rm ess}\lim_{x\to x_j^+}\int_0^{\tau}H_-(u_r(x,t)-k)\left[\varphi(u_r(x,t))-\varphi(k)\right ]\beta(t)\,dt\,\ge 0\,, \nonumber
\end{eqnarray}
since $\sgn(u)=1+2H_-(u)$. From inequalities \eqref{mkjn1} and \eqref{mkjn3} we obtain
$$
\begin{aligned}
&\int_0^{\tau}\!\!\!\!\int_{I_{j+1}} \!\! \big\{ [u_r-k]_+\zeta_t+ H_+( u_r-k) [\varphi(u_r)-\varphi(k)] \,\zeta_x\big\}\,dxdt +\\
+&\int_{I_{j+1}}[ u_{0r}(x)-k]_+\,\zeta(x,0)\,dx\,
 \ge   -\int_0^{\tau}\![h_j^+(t)-\varphi(k)] \,\zeta(x_j,t)\,dt\,.  
\end{aligned}
$$
Letting $k\to\infty$ in the above inequality gives
$$
\liminf_{k\to\infty}\int_0^{\tau}\![h_j^+(t)-\varphi(k)] \,\zeta(x_j,t)\,dt =
\int_0^{\tau}\![h_j^+(t)-\limsup_{k\to\infty}\varphi(k)] \,\zeta(x_j,t)\,dt \ge0\,,
$$
whence  the last inequality in \eqref{ext_h2} follows by the arbitrariness of $\zeta$.

Replacing $I_{j+1}\times(0,\tau)$ by $I_j\times(0,\tau)$, we obtain, similarly to  \eqref{mkjn1} and \eqref{mkjn3},
\begin{eqnarray}\label{mkjn2}
&&\ \int_0^{\tau}\!\!\!\int_{I_j}
 \big\{ (u_r-k)\zeta_t + [\varphi(u_r)-\varphi(k)] \,\zeta_x \big\}\,dxdt+ \int_{I_j}[ u_{0r}(x)-k]\,\zeta(x,0)\,dx\,=\\
&=& \!{\rm ess}\lim_{x\to x_j^-}\int_0^{\tau}\![\varphi(u_r(x,t))-\varphi(k)] \,\zeta(x,t)\,dt
= \int_0^{\tau}\![h_j^-(t)-\varphi(k)] \,\zeta(x_j,t)\,dt\,,  \nonumber
\end{eqnarray}
\begin{eqnarray}\label{mkjn4}
&&\ \int_0^{\tau}\!\!\!\int_{I_j} \left\{|u_r-k|\,\zeta_t+\sgn(u_r-k)\left [\varphi(u_r)-\varphi(k)\right ]\zeta_x\right\}dxdt \,+ \\
&+&\! \int_{I_j}|u_{0r}(x)-k|\,\zeta(x,0)\,dx\, - {\rm ess}\lim_{x\to x_j^-}\int_0^{\tau}\, \left[\varphi(u_r(x,t))-\varphi(k)\right ]\zeta(x,t)\,dt\,\ge 0\,, \nonumber 
\end{eqnarray}
whence 
$$
\begin{aligned}
&\int_0^{\tau}\!\!\!\int_{I_j}  \big\{ [u_r-k]_+\zeta_t+ H_+ (u_r-k) [\varphi(u_r)-\varphi(k)] \,\zeta_x\big\}\,dxdt\;+\\
 +\,&\int_{I_j}[ u_{0r}(x)-k]_+\,\zeta(x,0)\,dx\,\ge 
 \int_0^{\tau}\![h_j^-(t)-\varphi(k)] \,\zeta(x_j,t)\,dt  
\end{aligned}
$$
and $$
\limsup_{k\to\infty}\int_0^{\tau}\![h_j^-(t)-\varphi(k)] \,\zeta(x_j,t)\,dt =
\int_0^{\tau}\![h_j^-(t)-\liminf_{k\to\infty}\varphi(k)] \,\zeta(x_j,t)\,dt \le0\,.
$$
Since $\zeta$ is arbitrary we obtain the first inequality in \eqref{ext_h2}. 
\qed

\begin{remark} 
By standard density arguments and \eqref{lidesi1}, it follows from \eqref{mkjn4} that 
\begin{eqnarray}\label{nefo}
&&\int_0^{\tau}\!\!\!\int_{I_1} \left\{|u_r-k|\,\zeta_t+\sgn(u_r-k)\left [\varphi(u_r)-\varphi(k)\right ]\zeta_x\right\}dxdt \,+ \\
&+&\!\!\int_{I_1} |u_{0r}(x)-k|\,\zeta(x,0)\,dx \ge 
\int_0^{\tau}\! \left[h_1^-(t)-\varphi(k)\right ]\zeta(x_1,t)\,dt\nonumber 
\end{eqnarray}
if $\zeta\in C^1([0,\tau];C^1_c((-\infty,x_1]))$, $\zeta\ge0$, $\zeta(\cdot,\tau)=0$ in $(-\infty,x_1]$, and from \eqref{mkjn3} 
that 
\begin{eqnarray*}
&&\int_0^{\tau}\!\!\!\int_{I_{p+1}} \left\{|u_r-k|\,\zeta_t+\sgn(u_r-k)\left [\varphi(u_r)-\varphi(k)\right ]\zeta_x\right\}dxdt \,+ \\
&+&\!\! \int_{I_{p+1}}|u_{0r}(x)-k|\,\zeta(x,0)\,dx\, \ge -\int_0^{\tau}\! \left[h_p^+(t)-\varphi(k)\right ]\zeta(x_p,t)\,dt \nonumber 
\end{eqnarray*}
for all $\zeta\in C^1([0,\tau];C^1_c([x_p,\infty)))$, $\zeta\ge0$, $\zeta(\cdot,\tau)=0$ in $[x_p,\infty)$. Moreover, arguing as in the proof of Lemma \ref{incom} with $\eta_{\epsilon}$ in \eqref{defeta} replaced by
$$
\frac{2(x\!-\!x_j)\!-\!\epsilon}{\epsilon}\chi_{[x_j+\epsilon/2,x_j+\epsilon]}\,+\,\chi_{[x_j+\epsilon,x_{j+1}-\epsilon]}\,
+\frac{2(x_{j+1}\!-\!x)\!-\!\epsilon}{\epsilon}\chi_{[x_{j+1}-\epsilon,x_{j+1}-\epsilon/2]},
$$
we obtain that, for any $j=1,\dots,p-1$,
$$
\begin{aligned}
&\int_0^{\tau}\!\!\!\int_{I_{j+1}} \!\!\!\!\left\{|u_r\!-\!k|\,\zeta_t\!+\!\sgn(u_r\!-\!k)\left [\varphi(u_r)\!-\!\varphi(k)\right ]\zeta_x\right\}dxdt 
+\!\! \int_{I_{j+1}} \!\!\!\!|u_{0r}(x)\!-\!k|\,\zeta(x,0)\,dx\, \ge \\
\ge&\; -\!\int_0^{\tau}\! \left[h_j^+(t)-\varphi(k)\right ]\zeta(x_j,t)\,dt + \int_0^{\tau}\! \left[h_{j+1}^-(t)-\varphi(k)\right ]\zeta(x_{j+1},t)\,dt 
\end{aligned}
$$
for all $\zeta\in C^1([0,\tau];C^1_c([x_j,x_{j+1}]))$, $\zeta\ge0$, $\zeta(\cdot,\tau)=0$ in $[x_j,x_{j+1}]$. 
\end{remark} 
\begin{remark}\label{mer}  
For further reference we mention the following inequalities, which hold for all $\zeta\in C^1([0,\tau];C^1_c((-\infty,x_1]))$, $\zeta\ge0$, $\zeta(\cdot,\tau)=0$ in $(-\infty,x_1]$:
\begin{eqnarray}\label{subnefo}
&&\int_0^{\tau}\!\!\!\int_{I_1} \left\{[u_r-k]_+\,\zeta_t+H_+(u_r-k)\left [\varphi(u_r)-\varphi(k)\right ]\zeta_x\right\}dxdt \,+ \\
&+&\!\!\int_{I_1} [u_{0r}(x)-k]_+\,\zeta(x,0)\,dx\, \ge 
\int_0^{\tau}\! \left[h_1^-(t)-\varphi(k)\right ]\zeta(x_1,t)\,dt\,, \nonumber 
\end{eqnarray}
\begin{eqnarray}\label{supernefo}
&&\int_0^{\tau}\!\!\!\int_{I_1} \left\{[u_r-k]_-\,\zeta_t+H_-(u_r-k)\left [\varphi(u_r)-\varphi(k)\right ]\zeta_x\right\}dxdt \,+ \\
&+&\!\!\int_{I_1} [u_{0r}(x)-k]_-\,\zeta(x,0)\,dx\, \ge 0\,.\nonumber 
\end{eqnarray}
The proof is  analogous to that of \eqref{nefo}, starting from \eqref{subkru} and \eqref{superkru} instead of \eqref{mkru}. Similar inequalities hold in $S_j$ for $j=2,\dots,p+1$ (see Remark \ref{susu}).
\end{remark}


\subsection{Auxiliary problems}

Let $j=1,\dots,p+1$ and $n\in \N$. We consider the family of  auxiliary problems
$$
\begin{cases}
u_t+[\varphi(u)]_x=0&\text{in } S_j\\
u=n&\text{in } \partial I_j\times(0,T)\\
u=u_{0n}:=\min \{u_{0r},n\}&\text{in } I_j\times\{0\}\,. 
\end{cases} \leqno{(P_{j,n})}
$$
We follow \cite{MNRR,O} to define entropy solutions of $(P_{j,n})$.
\begin{definition}\label{otto}
 By an {\em entropy solution} of problem $(P_{j,n})$  we mean a function $u_{j,n}\in C([0,T]; L^1_{loc}(I_j))\cap  L^\infty(S_j)$ such that: 

\smallskip
\noindent $(i)$ $u_{j,n}$ is an entropy  solution of problem $(P)$ in $S_j$ (in the sense of Definition \ref{enso}) with Cauchy data $u_{0n}$;

\smallskip
\noindent $(ii)$ for all $\beta\in C^1_c(0,T)$, $\beta \geq 0$, $k\in[0,\infty)$ and $n\ge k$
\begin{equation}\label{dienzu}
\begin{aligned}
& {\rm ess} \!\!\!\lim_{x\to x_{j-1}^+}\! \int_0^T \!\!\! H_-(u_{j,n}(x,t)\!-\!k)[\varphi(u_{j,n}(x,t))\!-\!\varphi(k)]\beta(t)\,dt \le 0 \ \, \, \, \text{if $2\le j\le p\!+\!1$}\,,\\
& {\rm ess} \!\!\lim_{x\to x_j^-} \!\int_0^T \!\!\! H_-(u_{j,n}(x,t)-k)[\varphi(u_{j,n}(x,t))-\varphi(k)]\beta(t)\,dt \ge0
\ \, \, \, \text{if $1\le j\le p$}\,. 
\end{aligned}
\end{equation}
\end{definition}
\begin{remark}
By Definitions \ref{deso}-\ref{enso} and \ref{otto},  for all $\zeta\in C^1([0,T];C^1_c(I_j))$, $j,n$ and $k$ as above $u_{j,n}$ satisfies 
\begin{equation}\label{wefon}
\iint_{S_j} \big[u_{j,n}\zeta_t+ \varphi(u_{j,n}) \zeta_x\big]\,dxdt= - \int_{I_j}
u_{0n}(x)\zeta(x,0)\,dx\,,
\end{equation}
\begin{equation}\label{mkrun}
\iint_{S_j} \!\! \!\left\{|u_{j,n}\!-\!k|\,\zeta_t\!+\sgn(u_{j,n}\!-\!k)\left [\varphi(u_{j,n})\!-\!\varphi(k)\right ]\zeta_x\right\}dxdt 
\ge - \!\!\int_{I_j}\!\! |u_{0n}(x)\!-\!k|\,\zeta(x,0)dx. 
\end{equation}
By \eqref{superkru} and Remark \ref{susu}, there also holds
\begin{equation}\label{superkrun}
\iint_{S_j} \!\!\!\left\{[u_{j,n}\!-\!k]_-\,\zeta_t\!+\!H_-(u_{j,n}\!-\!k)\!\left [\varphi(u_{j,n})\!-\!\varphi(k)\right ]\zeta_x\right\}dxdt \ge -\!\! \!\int_{I_j}\!\![u_{0n}(x)-k]_-\,\zeta(x,0)dx 
\end{equation}
for all $\zeta\in C^1([0,\tau];C^1_c(\Omega))$, $\zeta\geq0$, $\zeta(\cdot,\tau)=0$ in $\Omega$, and $k\in[0,\infty)$. 
\end{remark}

\begin{prop}\label{pe1} 
Let $(H_0)$-$(H_1)$ hold. Then for all $j=1,\dots,p+1$ there exists an entropy solution  $u_j\in C([0,T]; L^1_{loc}(I_j))$ of problem $(P)$ in $S_j$, such that for all $k\ge 0$
\begin{equation}\label{dienzz}
\begin{aligned}
& {\rm ess} \lim_{x\to x_j^+} \int_0^T H_-(u_{j+1}(x,t)-k)[\varphi(u_{j+1}(x,t))-\varphi(k)]\beta(t)\,dt \le 0
\,, \\
& {\rm ess} \lim_{x\to x_j^-} \int_0^T H_-(u_j(x,t)-k)[\varphi(u_j(x,t))-\varphi(k)]\beta(t)\,dt \ge0\,. 
\end{aligned}
\end{equation}
Moreover, $u_j\in C([0,T]; L^1(I_j))$ for $j=2,\dots,p$.
\end{prop}
According to \cite{MNRR,O}, if $\varphi\in C^2([0,\infty))$ for every $j=1,\dots,p+1$, $n\in \N$ there exists a unique entropy solution $u_{j,n}$ of $(P_{j,n})$. To prove Proposition \ref{pe1} we need some preliminary results about these solutions. 
\begin{lem}\label{pe2}
Let $(H_1)$ hold, and let $u_{j,n}$ be the unique entropy solution of $(P_{j,n})$ $(j=1,\dots,p+1;n\in \N)$. Then:

\smallskip

\noindent $(i)$ there holds $0\le u_{j,n}\le n$, $u_{j,n}\le u_{j,n+1}$ a.e.~in $S_j$; 

\smallskip

\noindent $(ii)$ the sequence $\{u_{j,n}\}$ is bounded in $C([0,T]; L^1(I_j))$ if  $j=2,\dots,p$, $\{u_{1,n}\}$ is 
bounded in $C([0,T]; L^1_{loc}(I_1))$, and $\{u_{p+1,n}\}$ is bounded in $C([0,T]; L^1_{loc}(I_{p+1}))$.
\end{lem}
\proof
 $(i)$ We only give the proof if $j=1$. Consider the  problems
\begin{equation*}
\begin{cases}
u_t+[\varphi(u)]_x=0&\text{in }S_1\\
u=a_i&\text{in } 
\{x_1\}\times(0,T)\\
u=b_i &\text{in } 
I_1\times\{0\}\,,
\end{cases} \leqno{(P_i)}
\end{equation*}
where   $a_i \in L^{\infty}(0,T)$, $b_i\in L^{\infty}(I_1)$ $(i=1,2)$. As already mentioned, for each $i$ there exists  a unique entropy solution
 $z_i\in C([0,T]; L^1_{loc}(I_1))\cap  L^\infty(S_1)$  of $(P_i)$. Moreover,  for every $(x_0,t)\in S_1$ there holds (see \cite{MNRR, O})
\begin{equation}\label{ost}
\int_{x_0}^{x_1}\!\!|z_1(x,t)-z_2(x,t)|dx\le \!\int_{x_0-\|\varphi'\|_{\infty}t}^{x_1}\!\!|b_1(x)-b_2(x)|dx+
\|\varphi'\|_{\infty}\!\!\int_0^t\!\!|a_1(s)-a_2(s)|ds.
\end{equation}
Consider four sequences $\{a_{ik}\}\subset BV_{loc}(0,T)$, $\{b_{ik}\}\subset BV_{loc}(I_1)$ such that $a_{ik}\to a_i$ in $L^1_{loc}(0,T)$, $b_{ik}\to b_i$ in $L^1_{loc}(I_1)$ as $k\to\infty$ $(i=1,2)$. Let $z_{ik}\in BV_{loc}(S_1)$ be the unique entropy solution of $(P_i)$ with boundary and initial data $a_{ik}$, $b_{ik}$. As proven in \cite{Te}, for every $(x_0,t)\in S_1$ and $k\in\N$ there holds
\begin{eqnarray}\label{tc}
\int_{x_0}^{x_1}[z_{1k}(x,t)-z_{2k}(x,t)]_+\,dx&\le& \int_{x_0-\|\varphi'\|_{\infty}t}^{x_1}[b_{1k}(x)-b_{2k}(x)]_+\,dx\,+\\
&+&\|\varphi'\|_{\infty}\int_0^t[a_{1k}(s)-a_{2k}(s)]_+\,ds\,. \nonumber
\end{eqnarray}
On the other hand, applying \eqref{ost} to $z_{1k}$ and $z_{2k}\,$, by the arbitrariness of $x_0$ we obtain that  $z_{ik}(\cdot,t)\to z_i(\cdot,t)$ in $L^1_{loc}(I_1)$ as $k\to\infty$, for all $t\in(0,T)$. Hence there exists a subsequence $z_{i{k_l}}(\cdot,t)\subseteq z_{ik}(\cdot,t)$ such that $z_{i{k_l}}(\cdot,t)\to z_i(\cdot,t)$ a.e.~in $I_1$. Letting $k_l\to\infty$ in \eqref{tc} 
(with $k=k_l$), we obtain from  Fatou's Lemma that
\begin{eqnarray}\label{tca}
\int_{x_0}^{x_1}[z_1(x,t)-z_2(x,t)]_+\,dx&\le& \int_{x_0-\|\varphi'\|_{\infty}t}^{x_1}[b_1(x)-b_2(x)]_+\,dx\,+\\
&+&\|\varphi'\|_{\infty}\int_0^t[a_1(s)-a_2(s)]_+\,ds\,, \nonumber
\end{eqnarray}
whence the claim immediately follows.

\smallskip

\noindent  $(ii)$ Let $j\in\{2,\dots,p\}$ be fixed. Choosing in \eqref{wefon} $\zeta(x,t)=\alpha(x) \beta(t)$ with  $\alpha\in C^1_c(I_j)$, $\alpha\geq0$, and $\beta\in C^1([0,T])$, $\beta(T)=0$ we obtain
\begin{eqnarray*}
&& \left|\,  \int_0^T\!\!dt \,\beta'(t)\! \int_{I_j}\!\!
u_{j,n}(x,t)\,\alpha(x)\,dx\, \right |\le
\left| \, \int_0^T\!\!dt\, \beta(t) \! \int_{I_j}\!\!
\varphi(u_{j,n}(x,t)) \,\alpha'(x)\,dx\, \right | +\\
&+&\left|\, \beta(0) \int_{I_j}
u_{0n}(x)\,\alpha(x)\,dx\,\right|
\le \|\beta\|_{\infty}  \left(T \|\varphi\|_{\infty} + \|u_0\|_{L^1(I_j)} \right) \|\alpha\|_{W^{1,1}(I_j)}\,. \nonumber
\end{eqnarray*}
By standard smoothing arguments we can  set, for  fixed $\tau\in(0,T)$, $\beta=\beta_m$,
\begin{equation*}
\beta_m(t):= \chi_{(0,\tau]}(t)-m\left(t-\tau-\frac{1}{m}\right)\chi_{\left(\tau,\tau+\frac{1}{m}\right]}(t)
\quad \text{for }t\in [0,T]\,,
\end{equation*}
for  sufficiently large $m\in\N$. Letting $m\to\infty$ gives for all $\tau\in(0,T)$
\begin{equation}\label{zsa2}
0\le   \int_{I_j}
u_{j,n}(x,\tau)\,\alpha(x)\,dx\ \le \|\beta\|_{\infty}  \left(T \|\varphi\|_{\infty} + \|u_0\|_{L^1(I_j)} \right) \|\alpha\|_{W^{1,1}(I_j)}\,. 
\end{equation}
We fix $\epsilon>0$ and choose $\alpha$ in \eqref{zsa2} as 
$$
m\left(x\!-\!x_{j-1}\!-\! \epsilon\right) \chi_{\left(x_{j-1}+ \epsilon,x_{j-1}+ \epsilon+\frac{1}{m}\right]}
+ \chi_{(x_{j-1}+ \epsilon+\frac{1}{m},x_j - \epsilon-\frac{1}{m}]}
- m\left(x\!-\!x_j\!+ \!\epsilon\right)\chi_{\left(x_j - \epsilon-\frac{1}{m},x_j - \epsilon\right]}.
$$
Letting $m\to\infty$ we obtain that 
\begin{equation*}
0\le   \int_{x_{j-1}+ \epsilon}^{x_j- \epsilon} u_{j,n}(x,\tau)\,dx \le \|\beta\|_{\infty}  \left(T \|\varphi\|_{\infty} + \|u_0\|_{L^1(I_j)} \right) \left(x_j-x_{j-1}-2 \epsilon +2\right)\,, 
\end{equation*}
whence, by the arbitrariness of $\epsilon$,
\begin{equation}\label{zsa3}
0\le   \int_{I_j}
u_{j,n}(x,\tau)\,dx\ \le \|\beta\|_{\infty}  \left(T \|\varphi\|_{\infty} + \|u_0\|_{L^1(I_j)} \right) \left(x_j-x_{j-1}+2\right)\,. 
\end{equation}
This completes the proof if $j=2,\dots,p$. A similar argument can be used in  bounded subsets of $S_1$ and $S_{p+1}$, hence the conclusion follows.
\qed
\medskip

\noindent {\em Proof of Proposition \ref{pe1}.} By Lemma \ref{pe2}-$(i)$ we may define 
\begin{equation}\label{defuj}
u_j(x,t):=\lim_{n\to\infty} u_{j,n}(x,t) \quad \text{for a.e. $(x,t)\in S_j$}\,.
\end{equation}
Let $n\to\infty$ in \eqref{zsa3}. By monotonicity, $u_j\in L^{\infty}(0,T;L^1(I_j))$ and
\begin{equation}\label{c.nuo}u_{j,n}\to u_j\quad\mbox{in}\ \,L^1(S_j)
\end{equation}
for $j=2,\dots,p$. Similarly, 
$u_1\in L^{\infty}(0,T;L^1_{loc}(I_1))$, $u_{p+1}\in L^{\infty}(0,T;L^1_{loc}(I_{p+1}))$ and
\begin{equation}\label{c.nuo.bis}u_{1,n}\to u_1\ \ \mbox{in}\ \,L^1_{loc}(S_1),\quad u_{p+1,n}\to u_{p+1}\ \ \mbox{in}\ \,L^1_{loc}(S_{p+1})\,.
\end{equation}  
From the above convergences, letting $j\to\infty$ in \eqref{wefon} and \eqref{mkrun} we easily get
\begin{equation}\label{wefoj}
\iint_{S_j} \big[u_j\zeta_t+ \varphi(u_j) \zeta_x\big]\,dxdt= - \int_{I_j}
u_{0r}(x)\zeta(x,0)\,dx\,,
\end{equation}
\begin{equation}\label{mkruj}
\iint_{S_j} \left\{|u_j-k|\,\zeta_t+\sgn(u_j-k)\left [\varphi(u_j)-\varphi(k)\right ]\zeta_x\right\}dxdt 
\ge  - \int_{I_j}
|u_{0r}(x)-k|\,\zeta(x,0)\,dx
\end{equation}
 for all $\zeta\in C^1([0,\tau];C^1_c(I_j))$, $\zeta\geq0$, $\zeta(\cdot,\tau)=0$ in $I_j$. 

Next we show that  $u_j\in C([0,T]; L^1(I_j))$ for every $j=2,\dots,p$ (the same argument shows that $u_1 \in C([0,T]; L^1_{loc}(I_1))$ and $u_{p+1} \in C([0,T]; L^1_{loc}(I_{p+1}))$). By \cite[Proposition 3.10]{BSTT} and the above remarks there holds $u_j\in C((0,T]; L^1(I_j))$. To prove the continuity at $t=0$, observe that for any $\alpha\in C^1_c(I_j)$, $\alpha\geq 0$ and $h$ sufficiently small
\begin{eqnarray}\label{ec3}
&&\int_{I_j} |u_{j,n}(x,\tau)-u_{j,n}(x+h,\tau)|\,\alpha(x)\,dx \leq  \int_{I_j} |u_{0n}(x)-u_{0n}(x+h)|\,\alpha(x)\,dx + \\
&+&\!\! \int_0^{\tau}\!\!\!\int_{I_j} |\varphi(u_{j,n}(x,t))-\varphi(u_{j,n}(x+h,t))|\,|\alpha'(x)|\,dxdt \nonumber
\end{eqnarray}
for all $\tau\in (0,T)$ (the above inequality derives from the $L^1$-contraction property of the parabolic equation satisfied by the parabolic approximants of $u_{j,n}$; see \cite{MNRR, O}). 

By \eqref{c.nuo}-\eqref{c.nuo.bis}, as $n\to\infty$ in \eqref{ec3} we obtain for all $\tau\in (0,T)$ 
\begin{eqnarray}\label{ec4}
&&\int_{I_l} |u_j(x,\tau)-u_j(x+h,\tau)|\,\alpha(x)\,dx \leq \int_{I_l} |u_{0r}(x)-u_{0r}(x+h)|\,\alpha(x)\,dx + \\
&+&\!\!\int_0^{\tau}\!\!\!\int_{I_l}|\varphi(u_j(x,t))-\varphi(u_j(x+h,t))|\,|\alpha'(x)|\,dxdt\,. \nonumber
\end{eqnarray}
Let $\{\tau_k\}\subset(0,T)$, $\tau_k\to 0^+$ as $k\to\infty$. Since $u_{0r}\in L^1(I_j)$ and $\varphi(u_j)\in L^1(S_j)$, by \eqref{ec4} and the Fr\' echet-Kolmogorov Theorem  the sequence $\{u_j(\cdot,\tau_k)\,\alpha\}$ is relatively compact in $L^1(I_j)$. Then by \eqref{wefoj} and  a standard argument there holds $u_j(\cdot,\tau_k)\,\alpha  \to u_{0r}\,\alpha$ in $L^1(I_j)$ as $k\to\infty$\,. Arguing as in the proof of \cite[Proposition 3.10]{BSTT} we obtain that $\lim_{k\to \infty}  \int_{I_j}|u_j(x,\tau_k)-u_{0r}(x)|\,dx=0$, 
so $u_j\in C([0,T]; L^1(I_j))$.

\smallskip
It remains to prove \eqref{dienzz}. We only prove \eqref{dienzz}$_1$. By \eqref{defuj} and the Dominated Convergence Theorem, for a.e. $x\in I_j$ 
\begin{eqnarray*}
&&\lim_{n\to \infty}\int_0^T H_-(u_{j,n}(x,t)-k)[\varphi(u_{j,n}(x,t))-\varphi(k)]\beta(t)\,dt= \\
 &=&\int_0^T H_-(u_j(x,t)-k)[\varphi(u_j(x,t))-\varphi(k)]\beta(t)\,dt\,. \nonumber
\end{eqnarray*}
Then for a.e. $x\in I_j$ and every $\epsilon>0$ there exists $\bar n=\bar n(x)\ge k$ such that
\begin{eqnarray}\label{peq1}
&&\int_0^T H_-(u_j(x,t)-k)[\varphi(u_j(x,t))-\varphi(k)]\beta(t)\,dt \le \\
&\le&\!\!  \int_0^T H_-(u_{j,\bar n}(x,t)-k)[\varphi(u_{j,\bar n}(x,t))-\varphi(k)]\beta(t)\,dt+ \epsilon\,. \nonumber
\end{eqnarray}

On the other hand, arguing as in the proof of Lemma \ref{lim}, \eqref{superkrun} implies  that 
$$
x\mapsto \int_0^T\!\! H_-(u_{j,n}(x,t)-k)\left [\varphi(u_{j,n}(x,t))-\varphi(k)\right ]\beta(t)\,dt - kT \|\beta'\|_\infty\, x
$$
 is nonincreasing in $\R$.  Then by  \eqref{dienzu}$_1$ we get that for a.e. $x\in I_j$ and $n\ge k$
\begin{equation}\label{peq2}
\int_0^T\!\! H_-(u_{j,n}(x,t)-k)\left [\varphi(u_{j,n}(x,t))-\varphi(k)\right ]\beta(t)\,dt \le kT \|\beta'\|_\infty(x-x_{j-1})\,.
\end{equation}
By \eqref{peq1}-\eqref{peq2} and the arbitrariness of $\epsilon$, for a.e. $x\in I_j$ we obtain
\begin{equation*}
\int_0^T H_-(u_j(x,t)-k)[\varphi(u_j(x,t))-\varphi(k)]\beta(t)\,dt \le  kT \|\beta'\|_\infty(x-x_{j-1})\,,
\end{equation*}
whence \eqref{dienzz}$_1$ follows. 
\qed


\subsection{Existence proof}

Now we can prove Theorem \ref{exi}.

\smallskip

\noindent {\em Proof of Theorem \ref{exi}.}
Let  $u_j\in C([0,T]; L^1(I_j))$ $(2\le j \le p)$, $u_1\in C([0,T];L^1_{loc}(I_1))$ and $u_{p+1}\in C([0,T];L^1_{loc}(I_{p+1}))$ be given by Proposition \ref{pe1}, and let $h_j^{\pm}$ be given by Lemma \ref{incom}. For $j=1,\dots,p$ we set
\begin{equation}\label{C_j}
C_j(t):=\Big[\,  c_j - \int_0^t \left[h_j^+(s) - h_j^-(s) \right] ds
\,\Big]_+ \qquad (t\in [0,T])\,.
\end{equation}
Let $\bar{t}_j:=\sup \{\tau\in [0,T]\,|\, C_j(\tau)>0\}>0$. Then $\bar{t}_j>0$ since $C_j(0)= c_j>0$. By \eqref{ext_h2},  $C_j$ is nonincreasing in $(0,T)$, whence $C_j>0$ in $[0,\bar{t}_j)$ and, if $\bar{t}_j<T$, $C_j=0$ in $[\bar{t}_j,T]$ (observe that $\bar{t}_j=t_j$ for every $j=1,\dots,p$, with $t_j$ given by \eqref{dewait}). Let $\tau_1:=\min\,\{\bar{t}_1,\dots,\bar{t}_p \}$, and define $u\in C([0,\tau_1]; \ \mathcal M^+(\R))$ by setting
\begin{equation}\label{defuu}
\begin{cases}
u_r(\cdot,t):=u_j(\cdot,t) \text{\;\;in $I_j$} &(j=1,\dots, p+1) \\ 
u_s(\cdot,t):=\sum_{j=1}^p C_j(t)\delta_{x_j} 
\end{cases}
\qquad\text{for }0\le t\le \tau_1.
\end{equation}
It follows from Proposition \ref{pe1} that $u$ is an entropy solution of $(P)$ in $I_j\times(0,\tau_1)$ for $j=1,\dots,p+1$ which satisfies the compatibility condition  at every $x_1,\dots,x_p$ in $[0,\tau_1]$. Hence
$u$ is an entropy solution of $(P)$ in $\R\times(0,\tau_1)$, if we prove \eqref{ewf}-\eqref{mkru} with $\Omega=\R$, $\tau=\tau_1$ for all $\zeta\in C^1([0,\tau_1];C^1_c(\R))$, $\zeta\ge 0$, $\zeta(\cdot,\tau_1)=0$ in $\R$, such that 
$$
\text{${\rm supp}\,\zeta\cap \big( \{x_j\}\times (0,\tau_1)\big) \neq\emptyset$\;\; for some $j=1,\dots, p$\,.}
$$ 

We only give the proof when $\zeta(x,t)=\alpha(x) \beta(t)$ with  $\alpha\in C^1_c(\R)$, $\alpha\ge0$, $\alpha(x_j)>0$  for a unique $j\in\{1,\dots, p\}$, and $\beta\in C^1([0,\tau_1])$, $\beta\ge0$, $\beta(\tau_1)=0$ (the general case can be dealt with similarly). Let us first  prove \eqref{ewf} in this case, namely
\begin{eqnarray}\label{emo}
&&\int_0^{\tau_1}\!\!\!\int_{I_j\cup I_{j+1}}
 \big[u_r\zeta_t+ \varphi(u_r) \,\zeta_x\big]\,dxdt+ \int_{I_j\cup I_{j+1}} u_{0r}(x)\zeta(x,0)\,dx\,=\\
&-&\!\!\int_0^{\tau}\lla u_s(\cdot,t), \zeta_t(\cdot,t)\rra_{(x_{j-1},x_{j+1})}dt - \lla u_{0s},\zeta(\cdot,0)\rra_{(x_{j-1},x_{j+1})}\, \nonumber
\end{eqnarray}
for all $\zeta$ as above (we set $x_0:=-\infty$, $x_{p+1}:=\infty$).  From \eqref{C_j} and \eqref{defuu} we obtain \begin{eqnarray}\label{oz}
&&\int_0^{\tau_1}\lla u_s(\cdot,t), \zeta_t(\cdot,t)\rra_{(x_{j-1},x_{j+1})}dt+\lla u_{0s},\zeta(\cdot,0)\rra_{(x_{j-1},x_{j+1})}\,= \\
&=&\alpha(x_j)\left(\int_0^{\tau_1}\!\beta'(t)C_j(t)\,dt +c_j \beta(0)\! \right) 
=\alpha(x_j)\!\int_0^{\tau_1}\!\left[h_j^+(t) \!- \!h_j^-(t) \right]\beta(t)\,dt \nonumber
\end{eqnarray}
(see Lemma \ref{incom}). On the other hand, summing \eqref{mkjn1} and \eqref{mkjn2}  with $\tau=\tau_1$ gives
\begin{eqnarray}\label{emo1}
&&\int_0^{\tau_1}\!\!\!\int_{I_j\cup I_{j+1}}
 \big[u_r\zeta_t+ \varphi(u_r) \,\zeta_x\big]\,dxdt+ \int_{I_j\cup I_{j+1}} u_{0r}(x)\zeta(x,0)\,dx\,=\\
&=&-\,\alpha(x_j)\int_0^{\tau_1}\left[h_j^+(t) - h_j^-(t) \right]\beta(t)\,dt\,.\nonumber
\end{eqnarray}
and  \eqref{emo} follows from \eqref{oz} and \eqref{emo1}.

Next we  prove \eqref{mkru} for all $\zeta$ as above, namely 
\begin{eqnarray}\label{mkrum}
&&\int_0^{\tau_1}\!\!\!\int
_{I_j\cup I_{j+1}}
 \left\{|u_r-k|\,\zeta_t+\sgn(u_r-k)\left [\varphi(u_r)-\varphi(k)\right ]\zeta_x\right\}dxdt \,+\\
&+& \int_0^{\tau_1}\left\langle u_s(\cdot,t),\zeta_t(\cdot,t)\right\rangle_{(x_{j-1},x_{j+1})}\,dt +  \left\langle u_{0s}, \zeta(\cdot,0)\right\rangle_{(x_{j-1},x_{j+1})} \,\ge \nonumber\\
&\ge&-\int_{I_j\cup I_{j+1}}
|u_{0r}(x)-k|\,\zeta(x,0)\,dx \nonumber
\end{eqnarray}
Since  $u$ is  an entropy solution of $(P)$ in $I_j\times (0,\tau_1)$ and $I_{j+1}\times (0,\tau_1)$, and  satisfies the compatibility condition \eqref{coco} in $[0,\tau_1]$, it follows from \eqref{mkjn3} and \eqref{mkjn4}  that
\begin{eqnarray*}
&&\int_0^{\tau_1}\!\!\!\int
_{I_j\cup I_{j+1}}\!\! \left\{|u_r-k|\
\,\zeta_t+\sgn(u_r-k)\left [\varphi(u_r)-\varphi(k)\right ]\zeta_x\right\}dxdt \,+ \\
&+&\!  \int_{I_j\cup I_{j+1}}\!\!|u_{0r}(x)-k|\,\zeta(x,0)\,dx\, \ge\,
-\,\alpha(x_j)\int_0^{\tau_1}\left[h_j^+(t) - h_j^-(t) \right]\beta(t)\,dt\,. \nonumber
\end{eqnarray*}
Combined with \eqref{oz} this implies \eqref{mkrum}. Therefore, 
the measure $u$ defined by \eqref{defuu} is an entropy solution of $(P)$ in $\R\times(0,\tau_1)$. 

If $\tau_1<T$, either $u_s(\cdot,\tau_1)=0$, or $u_s(\cdot,\tau_1)>0$. If $u_s(\cdot,\tau_1)=0$, there holds $C_j(\tau_1)=0$ for all $j=1,\dots,p$ (see \eqref{C_j}-\eqref{defuu}), thus $ u_s(\cdot,t)=0$ for all $t\in[\tau_1,T]$. Then, by the standard theory of scalar conservation laws, we can continue the solution of $(P)$ in $(\tau_1,T]$ with initial data $u_r(\cdot,\tau_1)$. On the other hand,  if $ u_s(\cdot,\tau_1)>0$, then $C_j(\tau_1)>0$ for some $j=1,\dots,p$ and, arguing as before, we can continue the solution of $(P)$ in $(\tau_1,\tau_2]$, 
with initial data $u(\cdot,\tau_1)$,  for some $\tau_2\in (\tau_1,T]$. Iterating the procedure $q$ times with  $2\le q\le p$, we obtain that either $\tau_q=T$, or $u_s(\cdot,\tau_q)=0$. 
 \qed


\section{Proof of uniqueness}\label{unique} 
\setcounter{equation}{0}

This section is devoted to the proof of the uniqueness part of Theorem \ref{exiuni}: 

\begin{theorem}\label{uniq}
 Let $(H_0)$-$(H_1)$ be satisfied. Then there exists  at most one entropy solution of problem $(P)$, which belongs to  $C([0,T];\mathcal{M}^+(\R))$ and satisfies the compatibility condition at $x_j$ in $[0,t_j]$ for all $j=1,\dots,p$.
\end{theorem}
\proof Let $u, v\in C([0,T];\mathcal{M}^+(\Omega))$ be entropy solutions of $(P)$  satisfying the compatibility condition at every $x_j$ in $[0,t_j]$, and let
\begin{equation}\label{desta}
\tau:= \min\, \{t_u,t_v\} \quad \text{where}\quad
\begin{cases}
t_u:= \sup\, \{t\in[0,T) \,|\,  {\rm supp}\,u_s(\cdot,t)={\rm supp}\,u_{0s} \}\, \\
t_v:= \sup\, \{t\in[0,T) \,|\,  {\rm supp}\,v_s(\cdot,t)={\rm supp}\,u_{0s} \}\,.
\end{cases}
\end{equation}
Arguing as at the end of the proof of Theorem~\ref{exi}, it is enough to show that 
\begin{equation}\label{tp}
\text{ $u=v \quad$ in $\mathcal{M}(S_{\tau})$\,.}
\end{equation}
We claim that \eqref{tp} follows if we prove that 
\begin{equation}\label{contra}
\text{ $u_r=v_r$ \quad a.e.~in $\R\times(0,\tau)=:S_{\tau}$\,.}
\end{equation}
Indeed, \eqref{ewf} and \eqref{contra} imply that, for all $\zeta\in C^1([0,\tau];C^1_c(\R))$, $\zeta(\cdot,\tau)=0$ in $\R$,
$$
\int_0^{\tau}\!\! \langle u_s(\cdot,t) \! - \! v_s(\cdot,t), \zeta_t(\cdot,t)\rangle_{\R}\,dt \!
= \! \! \iint_{S_{\tau}} \! \{(u_r \! - \! v_r)\,\zeta_t+ [\varphi(u_r)\! - \! \varphi(u_r)] \zeta_x\}\,dxdt= 0
$$
Hence  $\langle u_s(\cdot,t)-v_s(\cdot,t),\alpha\rangle_{\R}=0$ for a.e.~$t\in (0,\tau)$, for all $\alpha \in C^1_c(\R)$.
Therefore $u_s=v_s$ in $L^{\infty}(0,\tau;\mathcal{M}(\R))$ and \eqref{tp} follows from \eqref{contra}.

It remains to prove \eqref{contra}, which is equivalent to showing that
\begin{equation}\label{contra1}
\text{ $u_r=v_r$ \quad a.e.~in $I_j\times(0,\tau)$ for all $j=1,\dots,p+1$\,.}
\end{equation}
We only prove \eqref{contra1} for $j=1$, since in the other cases the proof is similar.
Set $Q_1:=(-\infty,x_1]\times (0,\tau)$. We apply the Kru\v zkov method of doubling variables adapted to boundary valued problems (see \cite{MNRR, O, Se}). Let $\xi=\xi(x,t,y,s)$ be defined in $Q_1\times Q_1$, $\xi\ge0$, such that $\xi(\cdot,\cdot,y,s)\in C^1_c(Q_1)$ for every $(y,s)\in Q_1$ and  
$\xi(x,t,\cdot,\cdot)\in C^1_c(Q_1)$ for every $(x,t)\in Q_1$.
It follows from \eqref{nefo} that
\begin{eqnarray*}
&&\iint_{Q_1} \big\{\sgn(u_r(x,t)-v(y,s))[\varphi(u_r(x,t))-\varphi(v_r(y,s))]\xi_x(x,t,y,s) +  \\
&+&|u_r(x,t)\!-\!v_r(y,s)|\xi_t(x,t,y,s)\big\}\,dxdt\,\ge 
\int_0^{\tau}\! \!\left[h_1^-(t)\!-\!\varphi(v_r(y,s))\right ]\xi(x_1,t,y,s)dt\,, \nonumber 
\end{eqnarray*}
\begin{eqnarray*}
&&\iint_{Q_1} \big\{\sgn(u_r(x,t)-v(y,s))[\varphi(u_r(x,t))-\varphi(v_r(y,s))]\xi_y(x,t,y,s)+ \\
&+&|u_r(x,t)\!-\!v_r(y,s)|\xi_s(x,t,y,s)\big\}\, dyds \, \ge
\int_0^{\tau}\!\!\! \left[g_1^-(s)-\varphi(u_r(x,t))\right ]\xi(x,t,x_1,s)ds\,, \nonumber 
\end{eqnarray*}
where, by Lemma \ref{incom}-$(i)$, $g_j^{\pm}\in L^{\infty}(0,T)$ satisfies, for all $j=1,\dots,p$,  $g_j^{\pm}\ge0$ and
\begin{equation}\label{lidesi23}
{\rm ess} \lim_{x\to x_j^{\pm}}\int_0^T \varphi(v_r(x,t))\beta(t)\,dt=\int_0^T g_j^{\pm}(t)\beta(t)\,dt \quad\text{if }\beta\in L^1(0,T).
\end{equation}

Let $\rho_{\epsilon}$ $(\epsilon>0)$ be a symmetric mollifier in $\R$, and set in the previous inequalities 
\begin{equation}\label{ficho}
\xi(x,t,y,s)=\eta\,\Big(\frac{x+y}2,\frac{t+s}2\Big)\,\rho_{\epsilon}(x-y)\,\rho_{\epsilon}(t-s)
\end{equation} 
with $\eta\in C^1_c((-\infty,x_1]\times( 0,\tau))$, $\eta \ge0$. Then  we obtain
\begin{eqnarray}\label{xtys}
&&\iint\!\!\!\!\!\!\iint_{Q_1\times Q_1} \rho_{\epsilon}(x-y)\,\rho_{\epsilon}(t-s)
\Big\{|u_r(x,t)-v_r(y,s)| \, \eta_t\Big(\frac{x+y}2,\frac{t+s}2\Big) \, +\\
 &+&  \sgn(u_r(x,t)\!-\!v_r(y,s))[\varphi(u_r(x,t))\!-\!\varphi(v_r(y,s))]  \eta_x\Big(\frac{x\!+\!y}2,\frac{t\!+\!s}2\Big)\Big\} \,dxdtdyds \,\ge \nonumber \\
& \ge& \int_0^{\tau}\!\!\!\iint_{Q_1} \left[g_1^-(s)-\varphi(u_r(x,t))\right ] \,
\eta\,\Big(\frac{x+x_1}2,\frac{t+s}2\Big)\,\rho_{\epsilon}(x_1-x)\,\rho_{\epsilon}(t-s)\,dxdtds\,+ \nonumber\\
&+& \int_0^{\tau}\!\!\!\iint_{Q_1} \left[h_1^-(t)-\varphi(v_r(y,s))\right ]\,\eta\,\Big(\frac{x_1+y}2,\frac{t+s}2\Big)\,\rho_{\epsilon}(y-x_1)\,\rho_{\epsilon}(t-s)\,dydsdt\,. \nonumber
\end{eqnarray}
 
Concerning the right-hand side of \eqref{xtys}, by well-known properties of mollifiers 
$$
 \int_0^{\tau}\!\!\!\iint_{Q_1} g_1^-(s)\,\eta\,\Big(\frac{x+x_1}2,\frac{t+s}2\Big)\,\rho_{\epsilon}(x_1-x)\,\rho_{\epsilon}(t-s)\,dxdtds\, \to\, \frac 12\int_0^{\tau}\!\!g_1^-(s)\,\eta(x_1,s)\,ds\,,
$$
$$
 \int_0^{\tau}\!\!\!\iint_{Q_1} h_1^-(t)\,\eta\,\Big(\frac{x_1+y}2,\frac{t+s}2\Big)\,\rho_{\epsilon}(y-x_1)\,\rho_{\epsilon}(t-s)\,dydsdt\, \to\, \frac 12\int_0^{\tau}\!\!h_1^-(t)\,\eta(x_1,t)\,dt\,
$$
 as $\epsilon\to0^+$. Moreover, since 
$\iint_{Q_1}\rho_{\epsilon}(x_1-x)\,\rho_{\epsilon}(t-s)\,dxds =\tfrac 12$ for $\epsilon<\min\{t,\tau-t\}$,
$$
\begin{aligned}
&
\left | \int_0^{\tau}\!\!\!\!\iint_{Q_1}\!\!\!\varphi(u_r(x,t))\eta\Big(\frac{x\!+\!x_1}2,\frac{t\!+\!s}2\Big)
\rho_{\epsilon}(x_1\!-\!x)\rho_{\epsilon}(t\!-\!s)dxdtds \!-\! \frac 12\!\int_0^{\tau}\!\!\!h_1^-(t)\,\eta(x_1,t)dt \,\right |  \le\\
&\quad \le \iint_{Q_1}\!\! dxds\, \rho_{\epsilon}(x_1\!-\!x)\,\rho_{\epsilon}(t\!-\!s)  \int_0^{\tau}\!\!dt\,\varphi(u_r(x,t))\left |\,\eta\,\Big(\frac{x\!+\!x_1}2,\frac{t\!+\!s}2\Big)\!-\!\eta(x_1,t)\,\right | \,  +  \\
&\qquad + \left |\int_{I_1}dx\,\rho_{\epsilon}(x_1-x) \int_0^{\tau}[\varphi(u_r(x,t))-h_1^-(t)]\,\eta \,(x_1,t) 
 \,dt \int_0^{\tau}\rho_{\epsilon}(t-s)\,ds \,\right | 
 \le \\ 
&\quad \le  \frac{\|\varphi\|_{L^{\infty}}}{2} 
\sup_{0\le x_1-x+|t-s| \le \epsilon}\int_0^{\tau}\left|\,\eta\left(\frac{x_1+x}2,\frac {t+s} 2\right)-\eta(x_1,t)\right|\,dt+ \\
&\qquad + \|\rho_1\|_{\infty}\,\frac 1{\epsilon}\int_{x_1-\epsilon}^{x_1}\left|\int_0^{\tau}[\varphi(u_r(x,t))-h_1^-(t)]\,\eta\,(x_1,t)\,dt\,\right|\,dx\,. 
\end{aligned}
$$
By the smoothness of $\eta$ and equality \eqref{lidesi1}, the right-hand side of the above inequality vanishes as $\epsilon\to0^+$. Therefore, 
$$
\int_0^{\tau}\!\!\!\!\iint_{Q_1}\!\!\!\varphi(u_r(x,t))\eta\Big(\frac{x\!+\!x_1}2,\frac{t\!+\!s}2\Big)
\rho_{\epsilon}(x_1\!-\!x)\rho_{\epsilon}(t\!-\!s)dxdtds\,
\to\,\frac 12\int_0^{\tau}\!\!\!h_1^-(t)\,\eta(x_1,t)dt\, .
$$
It is similarly seen that
$$
\int_0^{\tau}\!\!\!\!\iint_{Q_1} \!\!\!\varphi(v_r(y,s))\eta\Big(\frac{x_1\!+\!y}2,\frac{t\!+\!s}2\Big)
\rho_{\epsilon}(y\!-\!x_1)\rho_{\epsilon}(t\!-\!s)dydsdt\,\to\,\frac12\int_0^{\tau}\!\!\!g_1^-(s)\eta\,(x_1,s)ds\,.
$$
Letting $\epsilon\to0^+$ in \eqref{xtys} we obtain that, for all $\eta\in C^1_c((-\infty,x_1]\times( 0,\tau))$, $\eta \ge0$,
\begin{eqnarray}\label{kruz}
&&\iint_{Q_1} \Big\{|u_r(x,t)-v_r(x,t)|\,\eta_t(x,t)+ \\
&+& \sgn(u_r(x,t)-v_r(x,t))[\varphi(u_r(x,t))-\varphi(v_r(x,t))]\eta_x(x,t)\Big\} dxdt\ge0\,. \nonumber
\end{eqnarray}

Now fix  $t'$, $t''$ such that  $0<t'<t''<\tau$, and $x_0<x_1$. Let  $\alpha_\delta=\alpha_\delta(x)$ and $\beta_\theta=\beta_\theta(t)$ be two families of mollifiers, such that $0<\delta<1$, $0<\theta<\min\{t',\tau-t''\}$. Set in  \eqref{kruz} 
$$
\eta(x,t)=\eta_{\delta,\theta}(x,t):=\int_{t-t''}^{t-t'} \beta_\theta(s)\,ds \int_{\|\varphi'\|_{\infty}
(t-t'')+x_0}^{x_1+\delta
}\!\! \alpha_\delta(x-y)\,dy \qquad(\delta>0)\,,
$$
with  $x\in (\|\varphi'\|_{\infty} (t-t'')+x_0,x_1]$, $t\in(t'-\theta,t''+\theta)$ (clearly, $\eta_{\delta,\theta}$ is nonnegative and belongs to $C^\infty_c((-\infty,x_1]\times(0,\tau)$)). Since $\alpha_\delta(x-x_1-\delta)=0$ if $x\in(-\infty,x_1]$, 
\begin{eqnarray*}
&&\iint_{Q_1}|u_r-v_r|\,[\beta_\theta(t-t')-\beta_\theta(t-t'')]\left(\int_{\|\varphi'\|_{\infty}
(t-t'')+x_0}^{x_1+\delta}\!\! \alpha_\delta(x-y)\,dy\right)\,dxdt\,-\\
&-&\!\!\iint_{Q_1}\{\|\varphi'\|_{\infty}
|u_r-v_r|+\sgn(u_r-v_r)[\varphi(u_r)-\varphi(v_r)]\}\,\times \nonumber\\
&\times&\!\!\! \alpha_\delta(x-\|\varphi'\|_{\infty}
(t-t'')-x_0)\left(\int_{t-t''}^{t-t'} \beta_\theta(s)\,ds\right)\,d xd t\ge0\,. \nonumber\\
\end{eqnarray*}
Since $ \|\varphi'\|_{\infty}
|u-v|+\sgn(u-v)[\varphi(u)-\varphi(v)]\ge0$  for every $u,v\ge0$, it follows that
\begin{equation*}
\iint_{Q_1}|u_r-v_r|\,[\beta_\theta(t-t')-\beta_\theta(t-t'')]
\left(\int_{\|\varphi'\|_{\infty}
(t-t'')+x_0}^{x_1+\delta
}\!\! \alpha_\delta(x-y)\,dy\right)dxdt\ge 0\,.
\end{equation*}

Let $\delta\to0^+$ in this inequality. Then, by the Dominated Convergence Theorem,  
$$
\int^{\tau}_0\int^{x_1}_{\|\varphi'\|_{\infty}
(t-t'')+x_0}|u_r(x,t)-v_r(x,t)|\,[\beta_\theta(t-t')-\beta_\theta(t-t'')])\,d xd t\ge 0\,,
$$
whence as $\theta\to0^+$
\begin{equation}\label{kruz2}
\int_{x_0}^{x_1}|u_r(x,t'')-v_r(x,t'')|\,dx\le\int^{x_1}_{x_0-\|\varphi'\|_{\infty}
(t''-t')}|u_r(x,t')-v_r(x,t')|\,dx.
\end{equation}
Since $u, v\in C([0,T];\mathcal{M}^+(\Omega))$, letting $t'\to 0^+$ in \eqref{kruz2} we obtain for all $(x_0,t'')\in Q_1$ 
\begin{equation*}
\int_{x_0}^{x_1} |u_r(x,t'')-v_r(x,t'')|\,dx\le \int_{x_0-\|\varphi'\|_{\infty}t''}^{x_1}|u_r(x,0)-v_r(x,0)|\,dx\,= 0\,,
\end{equation*}
since $u_r(\cdot,0)=v_r(\cdot,0)=u_{0r}$. Since $t''\in(0,\tau)$ is arbitrary, it follows that  $u_r= v_r$ in $Q_1$. 
\qed


\section{Comparison results}\label{comp} 
\setcounter{equation}{0}

To prove Theorem \ref{compa} we need  the following result.
\begin{lem}\label{curd}
Let $(H_1)$ be satisfied. Let $v_0\in\mathcal{M}^+(\R)$ satisfy $(H_0)$, and let $u_0\le v_0$ in $\mathcal{M}(\R)$. Let $u,v\in C([0,T];\mathcal{M}^+(\R))$ be the unique entropy solutions of $(P)$ with initial data $u_0, v_0$ given by Theorem \ref{exiuni}. Set
\begin{equation}\label{destabis}
\tau:= \min\, \{t_u,t_v\} \quad\text{where}\quad
\begin{cases}
t_u:= \sup\, \{t\in[0,T) \,|\,  {\rm supp}\,u_s(\cdot,t)={\rm supp}\,u_{0s} \}\,, \\
t_v:= \sup\, \{t\in[0,T) \,|\,  {\rm supp}\,v_s(\cdot,t)={\rm supp}\,v_{0s} \}\,.
\end{cases}
\end{equation}
Then 
\begin{equation}\label{disfobis}
\text{ $u_r(\cdot,t)\le v_r(\cdot,t)$ \; a.e. in $\R$  for any $t\in[0,\tau]$\,.}
\end{equation}
\end{lem}

\proof  
By assumption there holds 
$
u_{0s}=\sum_{j=1}^p c_j\delta_{x_j}$\,,  $v_{0s}=\sum_{j=1}^p d_j \delta_{x_j}$ 
with 
\begin{equation}\label{ius}
\text{$0\le  c_j\le d_j$\,, \; $d_j>0$ \;  for all $j=1,\dots,p$\,.} 
\end{equation}

Suppose first $c_j>0$ for all $j=1,\dots,p$.
 Let $u_{1,n}$ and $v_{1,n}$ be the entropy solutions of problem $(P_{1,n})$ with initial data $u_{0n}:=\min \{u_{0r},n\}$ and $v_{0n}:=\min \{v_{0r},n\}$, respectively. From inequality \eqref{tca} we get for any $(x_0,t)\in I_1\times[0,\tau)$
\begin{equation*}
\int_{x_0}^{x_1}[u_{1,n}(x,t)-v_{1,n}(x,t)]_+\,dx \le \int_{x_0-\|\varphi'\|_{\infty}t}^{x_1}[u_{0n}(x)-u_{0n}(x)]_+\,dx\,.
\end{equation*}
Since, by uniqueness and the proof of Theorem \ref{exi} (see \eqref{defuj} and \eqref{defuu}), 
$u_{1,n}\to u_r$ and $v_{1,n}\to v_r$  a.e. in $S_1$, we obtain from Fatou's Lemma that 
\begin{equation*}
\int_{x_0}^{x_1} [u_r(x,t)-v(x,t)]_+\,dx\le \int_{x_0-\|\varphi'\|_{\infty}t}^{x_1}[u_{0r}(x)-v_{0r}(x)]_+\,dx
\end{equation*}
for every $t\in(0,\tau)$. Similar inequalities can be proven in $I_j\times [0,\tau)$ for $j=2,\dots,p+1$, thus for every $x',x"\in \R$, $x'<x"$, and $t\in [0,\tau]$
\begin{equation*}
\int_{x'}^{x"} [u_r(x,t)-v_r(x,t)]_+\,dx\le \int_{x'-\|\varphi'\|_{\infty}t}^{x"+\|\varphi'\|_{\infty}t}[u_{0r}(x)-v_{0r}(x)]_+\,dx\,.
\end{equation*}
Hence the result follows  
in this case.

Now let $c_k=0$ for some $k\in\{1,\dots,p\}$; we only give the proof when $c_1=0$ and $c_j>0$ for $j=2,\dots,p$, since the general case is similar. Consider two sequences $\{u_{0m}\}, \{v_{0m}\}\subset BV_{loc}(\overline{I}_1\cup I_2)$ such that $u_{0m}\to u_{0r}$, $v_{0m}\to v_{0r}$ in $L^1_{loc}(\overline{I}_1\cup I_2)$ as $m\to\infty$, and $u_{0m}\le v_{0m}$ a.e. in $\overline{I}_1\cup I_2$ for all $m\in\N$. Set $u_{0m,n}:=\min\{u_{0m},n\}$, $v_{0m,n}:=\min\{v_{0m},n\}$ $(n\in\N)$, then consider the problems
\begin{equation*}
\begin{cases}
u_t+[\varphi(u)]_x=0&\text{in }(\overline{I}_1\cup I_2)\times(0,\tau)\\
u=n&\text{in } \partial(\overline{I}_1\cup I_2)\times(0,\tau)\\
u=u_{0m,n} &\text{in } (\overline{I}_1\cup I_2)\times\{0\}\,,
\end{cases} \leqno{(U)}
\end{equation*}
\begin{equation*}
\begin{cases}
v_t+[\varphi(v)]_x=0&\text{in }I_i\times(0,\tau)\\
v=n&\text{in } \partial I_i\times(0,\tau)\\
v=v_{0m,n} &\text{in } I_i\times\{0\}
\end{cases} \qquad (i=1,2)\,.
 \leqno{(V_i)}
\end{equation*}
Let $u_{m,n}\in BV_{loc}((\overline{I}_1\cup I_2)\times(0,\tau))$ and $v_{im,n}\in BV_{loc}(I_i\times(0,\tau))$ be the unique entropy solution of $(U)$ and $(V_i)$, respectively. Since $u_{m,n}\in BV_{loc}((\overline{I}_1\cup I_2)\times(0,\tau))$, for a.e. $t\in(0,\tau)$ there exist the traces $u_{m,n}(x_1^\pm ,t)$; moreover, there holds $u_{m,n}(x_1^\pm ,\cdot)\le n$ by comparison results (see the proof of Lemma \ref{pe2}). Then, since $u_{0m,n}\le v_{0m,n}$ a.e. in $\overline{I}_1\cup I_2$, it follows easily by comparison that 
\begin{equation}\label{comn}
\text{
$u_{m,n} \le v_{im,n}$ a.e. in $I_i\times(0,\tau)$ \qquad $(i=1,2)$\,.}
\end{equation}

As $m\to\infty$ there holds $u_{0m,n}\to u_{0n}$, 
$v_{0m,n}\to v_{0n}$ 
in $L^1_{loc}(\overline{I}_1\cup I_2)$. Then  
$u_{m,n}$ converges in $L^1_{loc}((\overline{I}_1\cup I_2)\times(0,\tau))$ to the unique entropy solution $u_n$ of the problem
\begin{equation*}
\begin{cases}
u_t+[\varphi(u)]_x=0&\text{in }(\overline{I}_1\cup I_2)\times(0,\tau)\\
u=n&\text{in } \partial(\overline{I}_1\cup I_2)\times(0,\tau)\\
u=u_{0n} &\text{in } (\overline{I}_1\cup I_2)\times\{0\}\,
\end{cases}
\end{equation*}
(see \cite{MNRR}).  Similarly, $v_{im,n}$ converges in $L^1_{loc}((I_i)\times(0,\tau))$ to the unique entropy solution $v_{in}$ of the problem
\begin{equation*}
\begin{cases}
v_t+[\varphi(v)]_x=0&\text{in }I_i\times(0,\tau)\\
v=n&\text{in } \partial I_i\times(0,\tau)\\
v=v_{0n} &\text{in } I_i\times\{0\}
\end{cases} \qquad (i=1,2)\,.
\end{equation*}
Then letting $m\to\infty$ in \eqref{comn} (possibly up to subsequences) we get
\begin{equation}\label{con}
\text{
$u_n \le v_{in}$ a.e. in $I_i\times(0,\tau)$ \qquad $(i=1,2)$\,.}
\end{equation}
By uniqueness and the proof of Theorem \ref{exi}, as $n\to\infty$ there holds $u_n\to u_r$ a.e. in $(\overline{I}_1\cup I_2)\times(0,\tau)$, $v_{in}\to v_r$  a.e. in $I_i\times(0,\tau)$ for $i=1,2$. Then from \eqref{con} letting $n\to\infty$ we obtain that $u_r \le v_r$ a.e. in $(\overline{I}_1\cup I_2)\times(0,\tau)$. It is similarly seen that $u_r \le v_r$ a.e. in $\bigcup_{j=3}^{p+1} I_j\times(0,\tau)$; since $u_r, v_r\in C([0,T];L^1_{loc}(\R))$, the result follows. 
\qed

\smallskip

\noindent {\em Proof of Theorem \ref{compa}.} By \eqref{disfobis} there holds 
\begin{equation}\label{disa}
\text{ $u_{ac}(\cdot,t)\le v_{ac}(\cdot,t)$ \;\; in $\mathcal{M}(\R)$ \;\; for any $t\in[0,\tau]$\,,}
\end{equation}
with $\tau$ given by \eqref{destabis}. Let us prove that
\begin{equation}\label{disb}
\text{ $u_s(\cdot,t)\le v_s(\cdot,t)$ \;\; in $\mathcal{M}(\R)$ \;\; for any $t\in[0,\tau]$\,.}
\end{equation}
We only prove \eqref{disb} when every $c_j$ is positive, since the proof is similar (and easier) if some $c_j$ is zero. 

To this purpose, as in the proof of Theorem \ref{uniq} we set $Q_1:=(-\infty,x_1]\times [0,\tau)$ and apply Kru\v zkov's method with $\xi=\xi(x,t,y,s)$ defined in $Q_1\times Q_1$, such that $\xi(\cdot,\cdot,y,s)\in C^1_c(Q_1)$ for every $(y,s)\in Q_1$ and  
$\xi(x,t,\cdot,\cdot)\in C^1_c(Q_1)$ for every $(x,t)\in Q_1$. From \eqref{subnefo}-\eqref{supernefo} we get
\begin{eqnarray*}
&&\iint_{Q_1} \big\{H_+(u_r(x,t)-v_r(y,s))[\varphi(u_r(x,t))-\varphi(v_r(y,s))]\,\xi_x(x,t,y,s)+\\
&&\ +[u_r(x,t)\!-\!v_r(y,s)]_+\xi_t(x,t,y,s)\big\}dxdt \ge \!\!
\int_0^{\tau}\! \!\!\left[h_1^-(t)\!-\!\varphi(v_r(y,s))\right ]\xi(x_1,t,y,s)dt, \nonumber 
\end{eqnarray*}
\begin{eqnarray*}
&&\iint_{Q_1} \big\{H_-(v_r(y,s)-u_r(x,t))[\varphi(v_r(y,s))-\varphi(u_r(x,t))]\,\xi_y(x,t,y,s)\,+ \\
&&\ +[v_r(y,s)\!-\!u_r(x,t)]_-\xi_s(x,t,y,s)\!\big\}dxdt =\!\! \iint_{Q_1}\!\!\!\! \big\{[u_r(x,t)\!-\!v_r(y,s)]_+\xi_s(x,t,y,s)+ \\
&&\ +H_+(u_r(x,t)-v_r(y,s))[\varphi(u_r(x,t))-\varphi(v_r(y,s))]\,\xi_y(x,t,y,s)\big\}\, dxdt \ge 0.\nonumber 
\end{eqnarray*}
Choosing $\xi$ as in \eqref{ficho}, from the above inequalities we get
\begin{eqnarray*}
&&\iint\!\!\!\!\!\!\iint_{Q_1\times Q_1}\rho_{\epsilon}(x-y)\,\rho_{\epsilon}(t-s) \Big\{[u_r(x,t)-v_r(y,s)]_+ \, \eta_t\Big(\frac{x+y}2,\frac{t+s}2\Big) 
 +\\
 &&\   +H_+(u_r(x,t)\!-\!v_r(y,s))[\varphi(u_r(x,t))\!-\!\varphi(v_r(y,s))] \, \eta_x\Big(\frac{x\!+\!y}2,\frac{t\!+\!s}2\Big)\Big\}dxdtdyds \! \ge
 \nonumber \\
&& \ \ge \int_0^{\tau}\!\!\!\!\iint_{Q_1} \!\!\!\left[h_1^-(t)\!-\!\varphi(v_r(y,s))\right ]\,\eta\,\Big(\frac{x_1\!+\!y}2,\frac{t\!+\!s}2\Big)\,\rho_{\epsilon}(y\!-\!x_1)\,\rho_{\epsilon}(t\!-\!s)\,dydsdt\,. \nonumber
\end{eqnarray*}
Then arguing as in the proof of Theorem \ref{uniq} plainly gives
\begin{eqnarray}\label{kruzz}
&&\iint_{Q_1} \Big\{H_+(u_r(x,t)-v_r(x,t))[\varphi(u_r(x,t))-\varphi(v_r(x,t))]\eta_x(x,t)+ \\
&&\ +[u_r(x,t)-v_r(x,t)]_+\,\eta_t(x,t)\Big\}\, dxdt
\ge \frac12\int_0^{\tau}[h_1^-(t)-g_1^-(t)]\,\eta(x_1,t)\,dt \nonumber
\end{eqnarray}
for every $\eta\in C^1_c((-\infty,x_1]\times( 0,\tau))$, $\eta \ge0$, with $h_1^-$, $g_1^-$ given by \eqref{lidesi} and \eqref{lidesi23}, respectively.

From \eqref{disfobis} and \eqref{kruzz} we obtain
\begin{equation}\label{disfoter}
\text{ $h_1^-\le g_1^-$ \;\; a.e in $(0,\tau)$.}
\end{equation}  
It is similarly proven that  (see Remark \ref{mer}):
\begin{equation}\label{cflu}
\text{$h_j^-\le g_j^-$ for $j=2,\dots,p$\,, \;\; $h_j^+\ge g_j^+$ for $j=1,\dots,p$ \quad  a.e. in $(0,\tau)$\,.}
\end{equation}

Now observe that, by uniqueness and the proof of Theorem \ref{exi},
$$
u_s(\cdot,t)=\sum_{j=1}^p C_j(t) \delta_{x_j}\,, \quad v_s(\cdot,t)=\sum_{j=1}^p D_j(t) \delta_{x_j}\,,
$$ 
with $C_j$ defined by \eqref{C_j} and
\begin{equation*} 
 D_j(t):=\left[d_j-\int_0^t [g_j^+(s)-g_j^-(s)]\,ds\right]_+ \qquad (j=1,\dots,p)\,.
\end{equation*}
By  \eqref{ius} and \eqref{cflu} there holds $C_j(t)\le D_j(t)$ for all $j=1,\dots,p$ and $t\in [0,\tau]$, thus inequality \eqref{disb} follows.

Hence there holds $u(\cdot,t)\le v(\cdot,t)$  in $\mathcal{M}(\R)$ for all $t\in[0,\tau]$. Arguing as at the end of the proof of Theorem~\ref{exi} we obtain the conclusion.
\qed

\begin{remark}
In section \ref{unique} we used Kru\v zkov's method to prove the uniqueness of entropy solutions satisfying the compatibility conditions. In the above proof  we used the same method to compare the fluxes of two such solutions at points where their singular parts are nontrivial and their regular parts are locally ordered. This additional information is contained in \eqref{disfoter}-\eqref{cflu}.
   \end{remark}




\begin{thebibliography}{999999}


\bibitem{AFP} 
\newblock L. Ambrosio, N. Fusco \& D. Pallara,
\newblock \emph{Functions of Bounded Variation and Free Discontinuity Problems},
\newblock (Oxford University Press, 2000).

\bibitem{BLN} 
C. Bardos, A. Y. Le Roux \& J. C. Nedelec,
 \emph{First order quasilinear equations with boundary condition},
 Comm. Partial Differential Equations \textbf{4} (1979), 1017-1034.

\bibitem{BSTT}
\newblock M. Bertsch, F. Smarrazzo, A. Terracina \& A. Tesei,
\newblock \emph{Radon measure-valued solutions of first order hyperbolic conservation laws},
\newblock preprint (2017).

\bibitem{CS}
\newblock G.-Q. Chen  \& Bo Su,
\newblock \emph{Discontinuous solutions of Hamilton-Jacobi equations: Existence, uniqueness, and
regularity,}
\newblock  Hyperbolic Problems: Theory, Numerics, Applications, T.Y. Hou et al. Eds. (Springer, 2003).


\bibitem{DS}
\newblock F. Demengel \& D. Serre,
\newblock \emph{Nonvanishing singular parts of measure valued solutions for scalar hyperbolic equations,}
\newblock Comm. Part. Diff. Equ. ~\textbf{16} (1991), 221-254.

\bibitem{Ev}
\newblock  L.C. Evans, 
\newblock \emph{Envelopes and nonconvex Hamilton-Jacobi equations,}
\newblock  Calc. Var \& PDE {\bf 50} (2014), 257-282. 

\bibitem{F}
\newblock	A. Friedman,
\newblock \emph{Mathematics in Industrial Problems, Part 8}, 
\newblock  IMA Volumes in Mathematics and its Applications {\bf 83} (Springer, 1997).

\bibitem{LP} 
\newblock  T.-P. Liu \& M.  Pierre, 
\newblock \emph{Source-solutions and asymptotic behavior
in conservation laws}, 
\newblock  J. Differential Equations \textbf{51} (1984), 419-441.

\bibitem{MNRR}
\newblock J. M\'alek, J. Ne\v cas, M. Rokyta \& M. R\. u\v zi\v cka, 
\newblock \emph{Weak and Measure-valued Solutions of Evolutionary PDEs}
\newblock (Chapman \& Hall, 1996).

\bibitem{O}
\newblock F. Otto, 
\newblock \emph{Initial-boundary value problem for a scalar conservation law}, 
\newblock Comptes Rendus Acad. Sci. \textbf{322} (1996), 729-734.

\bibitem{R1}
\newblock D.S. Ross, 
\newblock \emph{Two new moving boundary problems for scalar conservation laws}, 
\newblock Comm. Pure Appl. Math  \textbf{41} (1988), 725-737. 

\bibitem{R2}
\newblock D.S. Ross, 
\newblock \emph{Ion etching: An application of the mathematical theory of
hyperbolic conservation laws}, 
\newblock J. Electrochem. Soc. \textbf{135} (1988), 1235-1240.


\bibitem{Se} 
 \newblock  D. Serre,
\newblock \emph{ Systems of Conservation Laws, Vol. 1: Hyperbolicity, Entropies, Shock Waves},
\newblock (Cambridge University Press, 1999).

\bibitem{Te} A. Terracina,
 \emph{Comparison properties for scalar conservation laws with boundary conditions},
\newblock  Nonlinear Anal. \textbf{28} (1997), 633-
653.
 
\end{thebibliography}
\end{document}